\documentclass[10pt]{amsart}
\usepackage{indentfirst}
\usepackage[dvips]{graphicx}
\usepackage{graphicx}
\usepackage{newlfont}
\usepackage{amssymb}
\usepackage{amsmath,amscd}
\usepackage{latexsym}
\usepackage{amsthm}
\usepackage{psfrag}
\usepackage[usenames]{color}
\usepackage{textcomp}
\usepackage[all]{xy}

\newtheorem{thm}{Theorem}[section]

\theoremstyle{definition}
\newtheorem{defn}[thm]{Definition}

\theoremstyle{remark}

\newcommand{\R}{\mathbb R}
\newcommand{\N}{\mathbb N}
\newcommand{\Z}{\mathbb Z}

\newcommand{\p}{\varphi}
\newcommand{\s}{\psi}

\newcommand{\fr}{\vec \varphi}
\newcommand{\eps}{\varepsilon}

\numberwithin{equation}{section}

\begin{document}

\title[Multidimensional persistent homology and domain perturbations] {Stability of multidimensional persistent homology  with respect to domain perturbations}
\author{Patrizio Frosini}
\address{ Dipartimento di Matematica, Universit\`a di Bologna, P.zza di Porta S. Donato
5, I-$40126$ Bologna, Italia }
\email{frosini@dm.unibo.it}

\author{Claudia Landi}
\address{Dipartimento
di Scienze e Metodi dell'Ingegneria, Universit\`a di Modena e
Reggio Emilia, Via Amendola 2, Pad. Morselli, I-42100 Reggio
Emilia, Italia} \email{clandi@unimore.it}
\thanks{Research  partially carried out within the activities of ARCES (Universit\`a di Bologna).}


\date{} 



\begin{abstract}
Motivated by the problem of dealing with incomplete or imprecise  acquisition of data in computer vision and computer graphics,
 we extend results concerning the stability of persistent
homology with respect to function perturbations to results concerning the stability with respect to domain perturbations. Domain perturbations can be measured in a number of different ways. An important method to compare domains is the Hausdorff distance. We show that  by encoding sets using  the distance function,   the multidimensional matching distance between rank invariants of persistent homology groups is always upperly bounded by the Hausdorff distance between sets. Moreover, we prove that  our construction  maintains information about the original set.

Other well known methods to compare sets are considered, such as the symmetric difference distance between classical sets and the sup-distance between fuzzy sets. Also in these cases we present results stating that the multidimensional matching distance between rank invariants of persistent homology groups is upperly bounded by these distances.

An experiment showing the potential of our approach concludes the paper.
\end{abstract}

\maketitle

\section*{Introduction}
Persistent topology is a theory for studying objects related to computer vision and computer graphics, by analyzing the qualitative and quantitative behavior of real-valued functions defined over topological spaces and measuring the  shape properties of  the topological space under study (e.g., roundness, elongation, bumpiness, color). More precisely, persistent topology studies the sequence of nested lower level sets of the considered measuring functions and encodes  at which scale a topological feature  (e.g., a connected component, a tunnel, a void) is created, and when it is annihilated  along this filtration.  At the very beginning of the development of persistent topology, this encoding captured only the connected component changes in the lower level sets of  a real valued function, and took the name of {\em size function} \cite{Fr91,VeUrFrFe93}. Some years later, it was  extended to consider all homotopy groups of the lower level sets of a vector-valued function, under the name of {\em size homotopy groups} \cite{FrMu99}.  Nowadays we have a wide choice of variants for this encoding, ranging from persistent homology groups capturing the homology of a one-parameter increasing family of spaces \cite{EdLeZo02}, to  multidimensional persistent homology groups extending the previous concept to a multi-parameter setting \cite{CaZo09}, to vineyards coping with changes in the function over time \cite{CoEdMo06}, to interval persistence \cite{DeWe07}, just to cite a few. In this paper we focus on multidimensional persistent homology groups. For application purposes, these groups are further encoded by  considering only their rank, yielding  to a parametrized version of Betti numbers, called rank invariants (or persistent Betti numbers).

The stability of  multidimensional rank invariants is quite an important issue in persistent homology theory and its applications  because the lack of stability would make this invariant useless,  every data measurement being affected by noise. Stability with respect to perturbations of the measuring function  was proved in \cite{CeDi*09}, based on the results of \cite{BiCe*08,CaDiFe07}, comparing persistent homology groups by the multidimensional matching distance.

In this paper we consider the problem of stability with respect to changes of the topological space, which is as much important as the stability with respect to the change of measuring functions. Changes of the space under study can be measured in a number of different ways. Indeed, according to the kind of noise producing the perturbation, some distances are more suitable than other to compare sets. For example, the Hausdorff distance is useful to measure distortions of the domain, while the symmetric difference distance can cope with the presence of outliers. Due to the existence of many different ways to compare sets, we propose a general approach to the problem of stability of persistent homology groups with respect domain perturbations, and we apply this approach in a few cases.

Our main idea is to reduce the problem of stability with respect to changes of the topological space to that of stability with respect to changes of the measuring functions. This is achieved by  substituting the domain $K$ we are interested in with an appropriate function $f_K$ defined on a fixed set $D$ containing $K$, so that the perturbation of the set $K$ becomes a perturbation of the function $f_K$. As a consequence, the original measuring function $\vec \p_{|K}:K\rightarrow \R^k$ is replaced by a new measuring function $\vec \Phi:D\rightarrow \R^{k+1}$, $\vec\Phi=(f_K,\vec\p)$. Rank invariants of $(D,\vec\Phi)$ can be compared using the multidimensional matching distance.   In this way  we can prove  robustness of persistent homology groups under domain perturbations.

In particular, we use this strategy when sets are compared by the Hausdorff distance. In this case, taking $f_K$ equal to the distance function from $K$,  we prove that the multidimensional matching distance between the rank invariants associated with two compact sets $K_1$ and $K_2$ is always upperly bounded by the Hausdorff distance between $K_1$ and $K_2$ (Theorem \ref{stability_domains}). At the same time, we show that, in our approach, the information about the original domain $K$ and its original measuring function $\vec \p$ is fully maintained in the persistent homology groups of $(D,\vec\Phi)$ (Theorem \ref{lim}).

As a further contribution, we show stability with respect to perturbations of the domain measured using the symmetric difference distance. Also in this case, associating with $K$ a suitable function $f_K$, we can prove stability with  respect to domain perturbations measured by the symmetric difference distance (Theorem \ref{stability_domains_3}).

We also consider the situation where sets are described in a fuzzy sense, by means of probability density functions, easily obtaining a stability result also in this case (Theorem \ref{stability_domains_4}).

An experiment on a binary image concludes the paper, illustrating our results.

To conclude this introduction,  we emphasize three key-points in this paper. Firstly, there is a real need for developing techniques providing persistent homology with robustness against domain perturbations since the classical setting is not stable in this respect, as the example in Section \ref{example} clearly shows.
Secondly, the technique we propose in order to achieve stability makes sense only in the multidimensional version of persistent homology, since it is based on passing from a $k$-dimensional measuring function to a $(k+1)$-dimensional one. Thirdly, we underline that, substituting the study of the domain $K$ with that of a much simpler domain $D$, we obtain the searched for stability but not at the price of forgetting the persistent homology of $K$.   

\section{Preliminaries}\label{prel}

\subsection{Multidimensional persistent homology groups}

The following relations $\preceq$ and $\prec$ are defined in
$\R^k$: for $\vec u=(u_1,\dots,u_k)$ and $\vec v=(v_1,\dots,v_k)$,
we say $\vec u\preceq\vec v$ (resp. $\vec u\prec\vec v$) if and
only if $u_i\leq\ v_i$ (resp. $u_i<v_i$) for every index
$i=1,\dots,k$. Moreover, $\R^k$ is endowed with the usual
$\max$-norm: $\|(u_1,u_2,\dots,u_k)\|_{\infty}=\max_{1\leq i\leq
k}|u_i|$.

We shall use the following notations: $\Delta^+$ will be the open
set $\{(\vec u,\vec v)\in\R^k\times\R^k:\vec u\prec\vec v\}$.
Given a  topological space $X$, for every $k$-tuple $\vec
u=(u_1,\dots,u_k)\in\R^k$ and for every continuous function
$\vec\p:X\to\R^k$, we shall denote by $X\langle\fr\preceq \vec
u\,\rangle$ the {\em lower level  set} $\{x\in X:\varphi_i(x)\leq u_i,\
i=1,\dots,k\}$ and by $\|\vec \p\|_\infty$ the $\sup$-norm of $\vec \p$,
i.e. $\|\vec \p\|_\infty=\max_{x\in X}\|\vec \p(x)\|_\infty$.
The function $\vec\p$ will be called a $k$-dimensional \emph{measuring} (or {\em filtering}) {\em function}.

\begin{defn}
Let $\pi^{(\vec u,\vec v)}_q:\check{H}_q(X\langle\vec\p\preceq\vec
u\rangle)\rightarrow \check{H}_q(X\langle\vec\p\preceq\vec
v\rangle)$ be the homomorphism induced by the inclusion map
$\pi^{(\vec u,\vec v)}:X\langle\vec\p\preceq\vec
u\rangle\hookrightarrow X\langle\vec\p\preceq\vec v\rangle$ with
$\vec u\preceq\vec v$, where $\check{H}_q$ denotes the $q$th
\v{C}ech homology group. If $\vec u\prec\vec v$, the image of
$\pi^{(\vec u,\vec v)}_q$ is called the {\em multidimensional
$q$th persistent homology group of $(X,\vec\p)$ at $(\vec u, \vec
v)$}, and is denoted by $\check{H}_q^{(\vec u, \vec
v)}(X,\vec\p)$.
\end{defn}

In other words, the group $\check{H}_q^{(\vec u, \vec
v)}(X,\vec\p)$ contains all and only the homology classes of
$q$-cycles born before $\vec u$ and still alive at $\vec v$.

In what follows, we shall work with coefficients in a field
$\mathbb{K}$, so that homology groups are vector spaces, and hence
torsion-free. Therefore, they can be completely described by their
rank, leading to the following definition (cf. \cite{CaZo09}).

\begin{defn}[$q$th rank invariant]\label{Rank}
Let $X$ be a topological space, and $\fr:X\to\R^k$ a continuous
function. Let $q\in\mathbb{Z}$. The {\em $q$th rank invariant} of
the pair $(X,\fr)$  is the function
$\rho_{(X,\fr),q}:\Delta^+\to\N\cup\{\infty\}$ defined as
$$
\rho_{(X,\fr),q}(\vec u,\vec v)=\mathrm{rank}\,\mathrm{im}\,\pi^{(\vec u,\vec v)}_q.
$$
\end{defn}

 If $X$ is a triangulable space embedded in some $\R^n$, then  $\rho_{(X,\fr),q}(\vec u,\vec v)<+\infty$, for every $(\vec u,\vec v)\in\Delta^+$ and every $q\in\Z$ \cite{CaLaXX,CeDi*09}.

\subsection{Matching distance}\label{matchdist}

We now recall the construction of the distance $D_{match}$ to compare rank invariants of multidimensional persistent homology groups.

\begin{defn}[Admissible pairs]\label{Admissible}
\label{np} For every  vector $\vec{l}=(l_1,\ldots,l_k)$ of
$\mathbb{R}^k$ such that $l_i>0$ for $i=1,\dots,k$ and $\sum_{i=1}^kl_i^2=1$, and for every
vector $\vec{b}=(b_1,\ldots,b_k)$ of $\mathbb{R}^k$ such that
$\sum_{i=1}^k b_i=0$, we shall say that the pair
$(\vec{l},\vec{b})$ is \emph{admissible}. We shall denote the set
of all admissible pairs in $\R^k\times\R^k$ by $Adm_k$. Given an
admissible pair $(\vec{l},\vec{b})$, we define the half-plane
$\pi_{(\vec{l},\vec{b})}$ of $\R^k\times\R^k$ by the following
parametric equations:
$$
\left\{%
\begin{array}{ll}
    \vec u=s\vec l + \vec b\\
    \vec v=t\vec l + \vec b\\
\end{array}%
\right.
$$
for $s,t\in \R$, with $s<t$.
\end{defn}

The key property of this foliation of $\Delta^+$ by half-planes is  that the restriction of
$\rho_{(X,\fr),q}$ to each leaf can be seen as the
1-dimensional rank invariant associated with a suitable pair $(X,F_{(\vec l,\vec
b)}^{\fr})$, where $F_{(\vec
l,\vec b)}^{\fr}:X\rightarrow\R$.
Precisely, the following statement (proved in \cite{CaDiFe07}) holds:

\begin{thm}[Reduction Theorem]\label{Reduction}
Let $(\vec{l},\vec{b})$ be an admissible pair and let $F_{(\vec
l,\vec b)}^{\fr}:X\rightarrow\R$ be defined by setting
$$
F_{(\vec l,\vec
b)}^{\fr}(x)=\max_{i=1,\dots,k}\left\{\frac{\varphi_i(x)-b_i}{l_i}\right\}\
.
$$
Then, for every $(\vec u,\vec v)=(s\vec l+\vec b,t\vec l + \vec
b)\in\pi_{(\vec{l},\vec{b})}$ we have that
$$
\rho_{(X,\fr),q}(\vec u,\vec v)=\rho_{(X,F_{(\vec l,\vec
b)}^{\fr}),q}(s,t)\ .
$$
\end{thm}

Since 1-dimensional rank invariants can be compared by a distance $d_{match}$ that matches points of persistence diagrams \cite{CeDi*09,CoEdHa07}, we can construct the following distance between $k$-dimensional rank invariants, called {\em multidimensional matching distance}:

$$D_{match}(\rho_{(X,\fr),q},\rho_{(X,\vec\psi),q})=\sup_{(\vec
l,\vec b)\in Adm_k}\min_i l_i \cdot d_{match}\left(\rho_{(X,F_{(\vec
l,\vec b)}^{\fr}),q},\rho_{(X,G_{(\vec l,\vec
b)}^{\vec\psi}),q}\right), $$
where $F_{(\vec
l,\vec b)}^{\fr}=\max_{i=1,\dots,k}\left\{\frac{\varphi_i(x)-b_i}{l_i}\right\}$ and $G_{(\vec l,\vec
b)}^{\vec\psi}=\max_{i=1,\dots,k}\left\{\frac{\psi_i(x)-b_i}{l_i}\right\}$.

The key property of $D_{match}$ is its stability with respect to perturbations of the measuring function $\fr$, as the following theorem states (cf. \cite{CeDi*09}).

\begin{thm}[Multidimensional Stability Theorem]\label{Multidimensional}
Let $X$ be triangulable. For every $q\in\mathbb{Z}$, and for every two continuous functions $\fr,\vec\s:X\rightarrow \R^k$,
\begin{eqnarray*}
D_{match}\left(\rho_{(X,\fr),q},\rho_{(X,\vec\s),q}\right)\leq\|\fr-\vec\s\|_{\infty}.
\end{eqnarray*}
\end{thm}

\subsection{Comparison of sets}\label{distance}
The problems of description and comparison of sets can been dealt with in a myriad of different ways, each one more or less suitable than another for a given application task.

In classical set theory, the membership of elements in a set is assessed in binary terms according to a bivalent condition -- an element either belongs or does not belong to the set. By contrast, in fuzzy set theory \cite{Zadeh65, DuboisPrade80}, a fuzzy set $A$ in $X$ is characterized by a membership function $f_A:X\rightarrow [0,1]$, with the value $f_A(x)$ representing the grade of membership of $x$ in $A$. Usually, the nearer the value of $f_A(x)$ to $1$, the higher the grade of membership of $x$ in $A$.   The fuzzy set theory can be used in a wide range of domains in which information is incomplete or imprecise.

If classical set theory is adopted, then a number of different dissimilarity measures exist to compare two sets \cite{VeHa01,DeDe06}. A frequently used dissimilarity measure is the {\em Hausdorff distance}, which is defined for arbitrary non-empty compact subsets $K_1,K_2$ of  $\R^n$. If $K_1,K_2$ are contained in a non-empty compact subset $D$ of $\R^n$, the Hausdorff distance can be defined by
$$\delta_H(K_1,K_2)=\max\{\max_{x\in K_2}d_{K_1}(x), \max_{y\in K_1}d_{K_2}(y)\},$$
where $d_K$ denotes the distance to $K$, that  is  the function $d_K:D\to \R$ defined by $d_K(x)=\min_{y\in K}\|x-y\|$,  $\|\cdot\|$ being any norm on $\R^n$ (e.g., the Euclidean norm).  This can be reformulated as follows (cf. \cite[Ch.4, Sect. 2.2]{DeZo01}):
\begin{eqnarray}\delta_H(K_1,K_2)=\|d_{K_1}-d_{K_2}\|_\infty.\end{eqnarray}\label{hausinfty}

The Hausdorff distance is robust against small deformations, but it is sensitive to outliers: a single far-away noise point drastically increases the Hausdorff distance.

A dissimilarity measure that is based on the area of the symmetric difference, such as the symmetric difference pseudo-metric, overcomes the problem of outliers. Denoting by  $\mu$ the Lebesgue measure on $\R^n$, the {\em symmetric difference pseudo-metric} is defined between two measurable sets $A,B$ with finite measure by $d_\triangle(A,B)=\mu(A\triangle B)$ where $A\triangle B=(A\cup B)\setminus (A\cap B)$ is the symmetric difference of $A$ and $B$. It holds that $d_\triangle(A,B)=0$ if and only if $A$ and $B$ are equal almost everywhere. Identifying two sets $A$ and $B$ if $\mu(A\triangle B)=0$, we obtain the symmetric difference metric.

Other dissimilarity measures for more restricted patterns are, for example, the bottleneck distance between finite point sets and the Fr\'echet distance between curves.  However, since many other distances could be considered, we will limit our research to consider stability with respect to the Hausdorff and symmetric difference distances.

When fuzzy sets are used, their dissimilarity can be measured by any function distance. In this case we will confine  ourselves to consider the $\sup$-norm between fuzzy sets.

\section{Stability with respect to Hausdorff distance}\label{stabhaus}
Our main idea in proving stability of rank invariants with respect to noisy  domains  is to transform perturbations of sets into  perturbations of functions. In this way it is possible to apply the Multidimensional Stability Theorem \ref{Multidimensional}. When domain perturbations are measured by the Hausdorff distance, in order to pass from a set $K$  to a function, we insert the distance function $d_K$ described in subsection \ref{distance} as the first component of the measuring function. In this way,  assuming that all the sets under study are contained in a larger set $D$,  the original problem, i.e. studying persistent homology groups of
a set $K$ endowed with the restriction to $K$ of a measuring function $\vec \p:D\to \R^k$, is transformed into  the new problem of  studying the persistent homology groups of $D$ endowed with the measuring function $\vec\Phi=(d_K,\vec \p)$.

 Given two domains $K_1$ and $K_2$, and two functions $\vec\p_1,\vec\p_2:D\to \R^k$, our first result relates the distance $D_{match}$ between the new pairs $(D, \vec\Phi_1)$, $(D, \vec\Phi_2)$ to the change of the measuring functions $\vec \p_1$ and $\vec\p_2$, and to the Hausdorff distance between the original sets $K_1$, $K_2$. More precisely, it proves stability with respect to both set and function perturbations. Indeed, the change in  the multidimensional matching distance $D_{match}$ is shown to be never greater than the maximum among the change in the Hausdorff distance between the domains $K_1$ and $K_2$ and the change in the $\sup$-norm between the measuring functions $\vec \p_1$ and $\vec\p_2$. In particular, if $\vec \p_1$ and $\vec\p_2$ coincide then the change in  the multidimensional matching distance $D_{match}$ is  never greater than the Hausdorff distance between  $K_1$ and $K_2$.


\begin{thm}\label{stability_domains}
Let $K_1,K_2$ be non-empty closed subsets of a triangulable subspace $D$ of $\R^n$. Let $d_{K_1},d_{K_2}:D\to \R$ be their respective distance functions. Moreover, let $\vec \p_1,\vec\p_2 :D\to\R^k$ be  vector-valued continuous  functions. Then, defining $\vec \Phi_1,\vec \Phi_2:D\to \R^{k+1}$ by $\vec\Phi_1=(d_{K_1},\vec\p_1)$ and  $\vec\Phi_{2}=(d_{K_2},\vec\p_2)$, the following inequality holds:
$$D_{match}\left(\rho_{(D, \vec\Phi_1),q},\rho_{(D,\vec\Phi_2),q} \right)\le \max\left\{\delta_H(K_1,K_2),\|\vec\p_1-\vec\p_2\|_\infty\right\}.$$
\end{thm}

\begin{proof} The Multidimensional Stability Theorem \ref{Multidimensional} for measuring function perturbations implies that $D_{match}\left(\rho_{(D, \vec\Phi_1),q},\rho_{(D,\vec\Phi_2),q} \right)\le \|\vec\Phi_1-\vec\Phi_2\|_\infty$. It follows that $$D_{match}\left(\rho_{(D, \vec\Phi_1),q},\rho_{(D,\vec\Phi_2),q} \right)\le  \max\left\{\|d_{K_1}-d_{K_2}\|_\infty,\|\vec\p_1-\vec\p_2\|_\infty\right\}.$$  Hence, by equality (1.1), the claim is proved.
\end{proof}

We now consider the problem of retrieving  the rank invariants of $(K,\vec\p_{|K})$ from the rank invariants of $(D,\vec\Phi)$, with $\vec\Phi=(d_{K},\vec\p)$. The next result shows that for any sufficiently small value of $\beta\in\R$ there exists a sufficiently small value $\alpha\in\R$ with $0\le\alpha<\beta$ such that  $\rho_{(K,\vec \p_{|K}), q}(\vec u,\vec v)= \rho_{(D,\vec\Phi),q}\left((\alpha,\vec u),(\beta,\vec v)\right)$.

\begin{thm} \label{lim}
Let $K$ be  a non-empty  triangulable subset of a triangulable subspace $D$ of $\R^n$. Moreover, let $\vec\p:D\rightarrow \R^k$ be a continuous function.
Setting $\vec\Phi: D\to \R^{k+1}$, $\vec\Phi=(d_{K},\vec\p)$, for every $\vec u,\vec v\in \R^k$ with $\vec u\prec \vec v$, there exists a  real number $\hat\beta>0$ such that, for any $\beta\in \R$ with $0<\beta\le \hat\beta$, there exists a real number $\hat\alpha=\hat\alpha(\beta),$ with $0<\hat\alpha <\beta$, for which
$$\rho_{(K,\vec \p_{|K}), q}(\vec u,\vec v)= \rho_{(D,\vec\Phi),q}\left((\alpha,\vec u),(\beta,\vec v)\right),$$
for every $\alpha\in\R$ with $0\le \alpha\le \hat\alpha$. In particular,
\begin{eqnarray*}
\rho_{(K,\vec \p_{|K}), q}(\vec u,\vec v)= \lim_{\beta\to 0^+}\rho_{(D,\vec\Phi),q}\left((0,\vec u),(\beta,\vec v)\right).
\end{eqnarray*}\label{ptl}
\end{thm}

\begin{proof}
For every $\vec u\in\R^k$, we have
\begin{eqnarray*}
K\langle\vec\p_{|K}\preceq \vec u\rangle&=&\{x\in K:\vec\p(x)\preceq \vec u\}\\
&=&\{x\in D:d_K(x)\le 0\}\cap \{x\in D: \vec\p(x)\preceq \vec u\}\\
&=&\{x\in D:\vec\Phi(x)\preceq (0,\vec u)\}\\
&=& D\langle\vec\Phi\preceq (0,\vec u)\rangle.
\end{eqnarray*}
Hence, for every $q\in \Z$,   denoting by $\pi_q^{(\alpha,\vec u),(\beta,\vec v)}$ the homology homomorphism induced by the inclusion $D\langle\vec\Phi\preceq (\alpha,\vec u)\rangle\rightarrow D\langle\vec\Phi\preceq (\beta,\vec v)\rangle$, with $(\alpha,\vec u)\preceq (\beta,\vec v)$, it holds that
$$\rho_{(K,\vec\p_{|K}),q}(\vec u,\vec v)=\mathrm{rank\,}\left(\mathrm{im\,}\pi_q^{(0,\vec u),(0,\vec v)}\right).$$

We claim that there exists a positive real number $\hat\beta$ such that
$$\mathrm{im\,}\pi_q^{(0,\vec u),(0,\vec v)}\cong \mathrm{im\,}\pi_q^{(0,\vec u),(\beta,\vec v)}$$
 for every $\beta$ with $0<\beta\le \hat\beta$ (the claim is trivial for $\beta=0$). In particular, this fact proves that
$\rho_{(K,\vec\p_{|K}),q}(\vec u,\vec v)=\lim_{\beta\to  0^+}\rho_{(D,\vec\Phi),q}\left((0,\vec u),(\beta,\vec v)\right).$

In order to prove this claim, we consider the inverse system of homomorphisms $\pi_q^{(0,\vec u),(\beta,\vec v)}:\check{H}_q(D\langle\vec\Phi\preceq (0,\vec u\rangle)\to \check{H}_q(D\langle\vec\Phi\preceq (\beta,\vec v\rangle)$ over the directed set $\{\beta\in \R:\beta >0\}$ decreasingly ordered. The following  isomorphisms hold:
\begin{eqnarray*}
\mathrm{im\,}\pi_q^{(0,\vec u),(0,\vec v)}\cong  \mathrm{im\,}\varprojlim \pi_q^{(0,\vec u),(\beta,\vec v)}\cong \varprojlim\mathrm{im\,} \pi_q^{(0,\vec u),(\beta,\vec v)}.
\end{eqnarray*}
Indeed, $\mathrm{im\,}\pi_q^{(0,\vec u),(0,\vec v)}\cong  \mathrm{im\,}\varprojlim \pi_q^{(0,\vec u),(\beta,\vec v)}$ by the continuity of \v{C}ech homology, and $\mathrm{im\,}\varprojlim \pi_q^{(0,\vec u),(\beta,\vec v)}\cong \varprojlim\mathrm{im\,} \pi_q^{(0,\vec u),(\beta,\vec v)}$ because the inverse limit of vector spaces is an exact functor and therefore it preserves epimorphisms, and hence images.

It remains to prove that there exists a positive real number $\hat\beta$ such that, for every $0<\beta\le \hat \beta$, $\mathrm{im\,} \pi_q^{(0,\vec u),( \beta,\vec v)}$  is isomorphic to $\varprojlim\mathrm{im\,} \pi_q^{(0,\vec u),(\beta,\vec v)}$. To this end, let us consider the following commutative diagram, with $0<\beta'\le \beta''$:

\begin{eqnarray}\label{diagram}
\begin{array}{c}
\begin{centering}
\hfill \xymatrix {\check{H}_q(D\langle\vec\Phi\preceq (0,\vec u)\rangle)\ar[rr]^-{id} \ar[d]_{\pi_q^{(0,\vec u),(\beta',\vec v)}}&& \check{H}_q(D\langle\vec\Phi\preceq (0,\vec u)\rangle)\ar[d]^{\pi_q^{(0,\vec u),(\beta'',\vec v)}}\\
\check{H}_q(D\langle\vec\Phi\preceq (\beta',\vec v)\rangle)\ar[rr]^-{\pi_q^{(\beta',\vec v),(\beta'',\vec v)}}&&\check{H}_q(D\langle\vec\Phi\preceq (\beta'',\vec v)\rangle).}\hfill
\end{centering}
\end{array}
\end{eqnarray}
From the above diagram (\ref{diagram}), we see that each $\pi_q^{(\beta',\vec v),(\beta'',\vec v)}$ induces a map $\tau_q^{(\beta',\beta'')}:\mathrm{im\,}\pi_q^{(0,\vec u),(\beta',\vec v)}\rightarrow \mathrm{im\,}\pi_q^{(0,\vec u),(\beta'',\vec v)}$. From diagram (\ref{diagram}) we see that these maps are surjective. On the other hand, by the finiteness of the rank of  $\mathrm{im\,}\pi_q^{(0,\vec u),(0,\vec v)}$ and the monotonicity of the rank invariants, there exists $\hat \beta>0$ such that the rank of $\pi_q^{(0,\vec u),(\beta',\vec v)}$ is finite and equal to the rank of  $\pi_q^{(0,\vec u),(\beta'',\vec v)}$, whenever $0< \beta'\le \beta''\le\hat\beta$. Hence the maps  $\tau_q^{(\beta',\beta'')}$  are surjections between vector spaces of the same finite dimension, i.e. isomorphisms for every $0<\beta'\le \beta''\le \hat \beta$. Thus, $\varprojlim\mathrm{im\,} \pi_q^{(0,\vec u),(\beta,\vec v)}$ is the inverse limit of a system of  finite dimensional vector spaces isomorphic to  $\mathrm{im\,} \pi_q^{(0,\vec u),(\hat\beta,\vec v)}$, proving
 the claim.

We now claim that for every strictly positive real number $\beta$,  there exists a strictly positive real number $\hat\alpha<\beta$ such that
$$\mathrm{im\,}\pi_q^{(0,\vec u),(\beta,\vec v)}\cong \mathrm{im\,}\pi_q^{(\alpha,\vec u),(\beta,\vec v)}$$
 for every $\alpha$ with $0\le \alpha\le \hat\alpha$.

 This claim can be proved in much the same way as the previous one. We consider the inverse system of homomorphisms $\pi_q^{(\alpha,\vec u),(\beta,\vec v)}:\check{H}_q(D\langle\vec\Phi\preceq (\alpha,\vec u\rangle)\to \check{H}_q(D\langle\vec\Phi\preceq (\beta,\vec v\rangle)$ over the directed set $\{\alpha\in \R:0\le\alpha<\beta\}$ decreasingly ordered. The following  isomorphisms follow again from  the continuity of \v{C}ech homology and the exacteness of  the inverse limit functor for  vector spaces:
\begin{eqnarray*}
\mathrm{im\,}\pi_q^{(0,\vec u),(\beta,\vec v)}\cong  \mathrm{im\,}\varprojlim \pi_q^{(\alpha,\vec u),(\beta,\vec v)}\cong \varprojlim\mathrm{im\,} \pi_q^{(\alpha,\vec u),(\beta,\vec v)}.
\end{eqnarray*}

To prove that there exists a strictly positive real number $\hat\alpha$ such that, for every $0\le\alpha\le \hat\alpha$, $ \mathrm{im\,} \pi_q^{(\alpha,\vec u),( \beta,\vec v)}$  is isomorphic to $\varprojlim\mathrm{im\,} \pi_q^{(\alpha,\vec u),(\beta,\vec v)}$, let us consider the following commutative diagram, with $0\le \alpha'\le \alpha''$:

\begin{eqnarray}\label{diagram2}
\begin{array}{c}
\begin{centering}
\hfill \xymatrix {\check{H}_q(D\langle\vec\Phi\preceq (\alpha',\vec u)\rangle)\ar[rr]^-{\pi_q^{(\alpha',\vec u),(\alpha'',\vec u)}} \ar[d]_{\pi_q^{(\alpha',\vec u),(\beta,\vec v)}}&& \check{H}_q(D\langle\vec\Phi\preceq (\alpha'',\vec u)\rangle)\ar[d]^{\pi_q^{(\alpha'',\vec u),(\beta,\vec v)}}\\
\check{H}_q(D\langle\vec\Phi\preceq (\beta,\vec v)\rangle)\ar[rr]^-{id}&&\check{H}_q(D\langle\vec\Phi\preceq (\beta,\vec v)\rangle).}\hfill
\end{centering}
\end{array}
\end{eqnarray}
From the above diagram (\ref{diagram2}), we see that each $\pi_q^{(\alpha',\vec u),(\alpha'',\vec u)}$ induces a map $\sigma_q^{(\alpha',\alpha'')}:\mathrm{im\,}\pi_q^{(\alpha',\vec u),(\beta,\vec v)}\rightarrow \mathrm{im\,}\pi_q^{(\alpha'',\vec u),(\beta,\vec v)}$. From diagram (\ref{diagram2}) we see that these maps are injective. On the other hand, by the finiteness of the rank of $\mathrm{im\,}\pi_q^{(\alpha,\vec u),(\beta,\vec v)}$, for any $\alpha$ with $0<\alpha<\beta$, and the monotonicity of the rank invariants, there exists $\hat \alpha$, with $0< \hat\alpha<\beta$, such that the rank of $\pi_q^{(\alpha',\vec u),(\beta,\vec v)}$ is finite and equal to the rank of  $\pi_q^{(\alpha'',\vec u),(\beta,\vec v)}$, whenever $0\le  \alpha'\le \alpha''\le\hat\alpha$. Hence the maps  $\sigma_q^{(\alpha',\alpha'')}$  are injections between vector spaces of the same finite dimension, i.e. isomorphisms for every $0\le \alpha'\le \alpha''\le \hat \alpha$. Thus, $\varprojlim\mathrm{im\,} \pi_q^{(\alpha,\vec u),(\beta,\vec v)}$ is the inverse limit of a system of  finite dimensional vector spaces isomorphic to  $\mathrm{im\,} \pi_q^{(\hat\alpha,\vec u),(\hat\beta,\vec v)}$, proving
 the claim.
\end{proof}

Many applications require that the presence of single outliers does not affect the evaluation of similarity.  In these cases, always assuming $K$ triangulable,  it is sufficient to  study the closure of the interior of $K$ instead of $K$ itself. Indeed, applying Theorems \ref{stability_domains} and \ref{lim}  with the closure of the interior of $K$ instead of $K$, we obtain a result of stability of persistent homology groups with respect to the perturbations of the studied set and a reconstruction result for the original   persistent homology groups modulo perturbations of zero measure.

We underline once more that the results of this section are based  on the idea of translating the problem of stability with respect to set perturbations into that of stability with respect to function perturbations. Therefore, the use of the distance function  is only one among many ways to achieve this end and has the advantage of working well when sets are compared using the Hausdorff distance. One could conceive different ways, in connection with other methods to compare sets, as the following sections show.

\section{Stability with respect to other distances between sets}\label{symmetric}
Our approach can be easily adapted to noise that is small with respect to distances other than the Hausdorff distance $\delta_H$.

We first show how persistent homology  can be made stable with respect to perturbations of the sets measured using the symmetric difference distance (Theorem \ref{stability_domains_3}). Then we show the stability with respect to perturbations of fuzzy sets (Theorem \ref{stability_domains_4}).

\subsection{Stability with respect to the symmetric difference distance}
We work with a non-empty closed subset $K$ of  a  triangulable set $D$ in $\R^n$.
In this case, instead of the distance function $d_K$, our construction depends on the use of functions
$\lambda_K^\eps:\R^n\to \R$, with $\eps\in\R$, $\eps>0$, defined as
$$\lambda_K^\eps(x)=\mu(B_\eps)^{-1}\cdot \int_{y\in B_\eps(x)}\chi_K(y)\ \mathrm{d}y$$
where $B_\eps(x)$ denotes the $n$-disk centered at $x$ with radius $\eps$, $B_\eps=B_\eps(\vec 0)$,   and $\chi_K$ denotes the characteristic function of $K$.  The underlying idea of this choice is that  the closer a point $x$ of the real plane is to a large part of $K$, the closer the value of $\lambda_K^\eps(x)$ is to $1$. More precisely, $\lambda_K^\eps(x)=1$ if and only if $\mu(B_\eps(x)\cap K)=\mu(B_\eps)$, whereas $\lambda_K^\eps(x)=0$ if and only if $\mu(B_\eps(x)\cap K)=0$. Clearly, $\lambda_K^\eps$ is a continuous function for every real number $\eps>0$.

Analogously to Theorem \ref{stability_domains}, in this case we have the following result.

\begin{thm}\label{stability_domains_3}
Let $K_1,K_2$ be non-empty closed subsets of a  triangulable subspace $D$ of $\R^n$. Moreover, let $\vec \p_1,\vec\p_2 :D\to\R^k$ be  vector-valued continuous  functions. Then, defining $\vec \Psi_1^{\eps},\vec \Psi_2^{\eps}:D\to \R^{k+1}$ by $\vec\Psi_1^\eps=(-\lambda^\eps_{K_1},\vec\p_1)$ and  $\vec\Psi_{2}^\eps=(-\lambda^\eps_{K_2},\vec\p_2)$, the following inequality holds:
\begin{eqnarray}
D_{match}\left(\rho_{(D, \vec\Psi_1^\eps),q},\rho_{(D,\vec\Psi_2^\eps),q} \right)\le \max\left\{ \frac{d_\triangle(K_1,K_2)}{\mu(B_\eps)},\|\vec\p_1-\vec\p_2\|_\infty\right\}.
\end{eqnarray}\label{inequality}
\end{thm}

\begin{proof}
For every $x\in D$,
\begin{eqnarray*}
|\lambda^\eps_{K_1}( x)-\lambda^\eps_{K_2}( x)|
& =& \mu(B_\eps)^{-1}\cdot \left|\int_{y\in B_\eps(x)} \chi_{K_1}( y)-\chi_{K_2}( y)\ \mathrm{d}y\right| \\
  &\le& \mu(B_\eps)^{-1}\cdot \int_{D} |\chi_{K_1}( y)-\chi_{K_2}( y)|\ \mathrm{d}y \\
  &=&  \mu(B_\eps)^{-1}\cdot \mu(K_1\triangle K_2).
\end{eqnarray*}
Thus $\| \lambda^\eps_{K_1}-\lambda^\eps_{K_2}\|_\infty\le   \mu(B_\eps)^{-1}\cdot \mu(K_1\triangle K_2)$.
The Multidimensional Stability Theorem \ref{Multidimensional} for measuring function perturbations implies that
$$D_{match}\left(\rho_{(D, \vec\Psi_1^\eps),q},\rho_{(D,\vec\Psi_2^\eps),q} \right)\le \|\vec\Psi_1^\eps-\vec\Psi_2^\eps\|_\infty.$$ It follows that
\begin{eqnarray*}
D_{match}\left(\rho_{(D, \vec\Psi_1^\eps),q},\rho_{(D,\vec\Psi_2^\eps),q} \right)&\le&  \max\left\{\|\lambda_{K_1}^\eps-\lambda_{K_2}^\eps\|_\infty,\|\vec\p_1-\vec\p_2\|_\infty\right\}\\
&\le& \max\left\{\mu(B_\eps)^{-1}\cdot \mu(K_1\triangle K_2),\|\vec\p_1-\vec\p_2\|_\infty\right\}\\
& = & \max\left\{\mu(B_\eps)^{-1}\cdot d_\triangle(K_1,K_2),\|\vec\p_1-\vec\p_2\|_\infty\right\}.
\end{eqnarray*}
\end{proof}

The previous theorem shows that, under our hypotheses, if two compact subsets $K_1,K_2$ of the real plane are close to each other in the sense that their symmetric difference has a small measure, then also the rank invariants constructed using the functions $\vec\Psi_{1}^\eps$, $\vec\Psi_{2}^\eps$ are close to each other.

We observe that the estimate in  inequality (\ref{inequality}) can be improved by substituting  $d_\triangle(K_1,K_2)$ with $\max_{x\in \R^n}\left|\int_{y\in B_\eps(x)} \chi_{K_1}( y)-\chi_{K_2}( y)\ \mathrm{d}y\right|$.

\subsection{Stability with respect to perturbations of fuzzy sets}

Now we consider  the case when sets are defined according to  fuzzy theory, that is  through functions representing the grade of membership of points to the considered set. One obtains a fuzzy set, for example, when a probability density $p(x)$ is given, $p(x)$ expressing the probability that a point of the considered set belongs to an infinitesimal neighborhood of $x$. We confine ourselves to considering only probability densities with compact support contained in a triangulable subspace $D$ of $\R^n$. From the Multidimensional Stability Theorem \ref{Multidimensional} for measuring function perturbations we immediately deduce the following result, whose simple proof is omitted, concerning the stability with respect to perturbations of fuzzy sets defined by probability densities.

  \begin{thm}\label{stability_domains_4}
Let $p_1,p_2$ be two probability density functions having support contained in a compact and triangulable subspace $D$ of $\R^n$. Defining $\vec \Psi_1,\vec \Psi_2:D\to \R^{k+1}$ by $\vec\Psi_1=(-p_1,\vec\p_1)$ and  $\vec\Psi_{2}=(-p_2,\vec\p_2)$, the following statement holds:
$$D_{match}\left(\rho_{(D, \vec\Psi_1),q},\rho_{(D,\vec\Psi_2),q} \right)\le \max\left\{\|p_1-p_2\|_\infty,\|\vec\p_1-\vec\p_2\|_\infty\right\}.$$
\end{thm}

\section{An example}\label{example}

In this section the theoretical framework presented in   Section \ref{stabhaus} is applied in a discrete setting. Our goal is to
check the stability of the proposed framework with respect to set perturbations measured by the Hausdorff distance. We confine ourselves to the case $q=0$.

\begin{figure}[h]
\begin{center}
\begin{tabular}{cc}\fbox{\includegraphics[width=5cm]{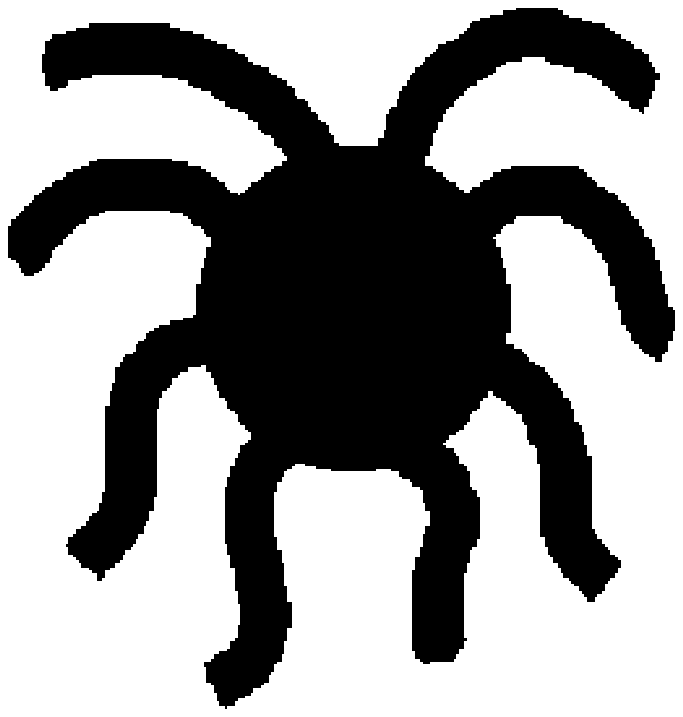}} &  \fbox{\includegraphics[width=5cm]{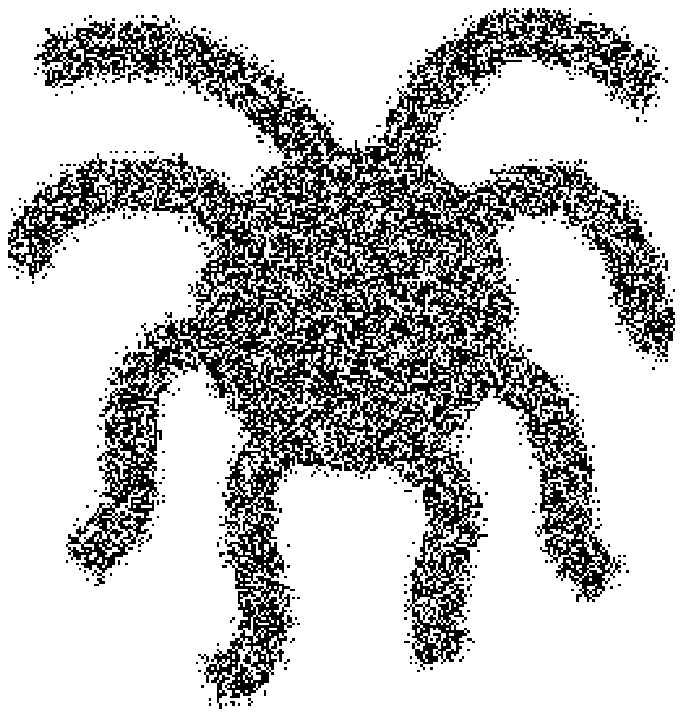}}\end{tabular}
\caption{Two binary images of an octopus. The image on the right is a noisy version of that on the left.}\label{octopus}
\end{center}
\end{figure}

With this aim in mind, we work with the binary  digital image represented in Figure \ref{octopus} (left), and we corrupt this image by adding {\em salt \& pepper} noise to a neighborhood of the set of its black pixels, as shown in Figure \ref{octopus} (right).

Black pixels of left and right images represent the sets $K_1$, $K_2$ under study, respectively, whereas in both cases the 269x256 rectangle of black and white pixels together constitute the set $D$. The so obtained noisy set $K_2$ is close to the original set $K_1$ with respect to the Hausdorff distance.

A graph structure
based on the local $4$-neighbors adjacency relations of the digital points is used in order to topologize  the images.

Fixed the point $c\in D$ corresponding to the center of mass of $K_1$, the chosen measuring function for both instances is $\p:D\to \R$, $\p(p)=-\|p-c\|$. 

\begin{center}
\begin{figure}[h]
\psfrag{u}{$u$}
\psfrag{v}{$v$}
\begin{tabular}{ll}\includegraphics[width=0.5\textwidth]{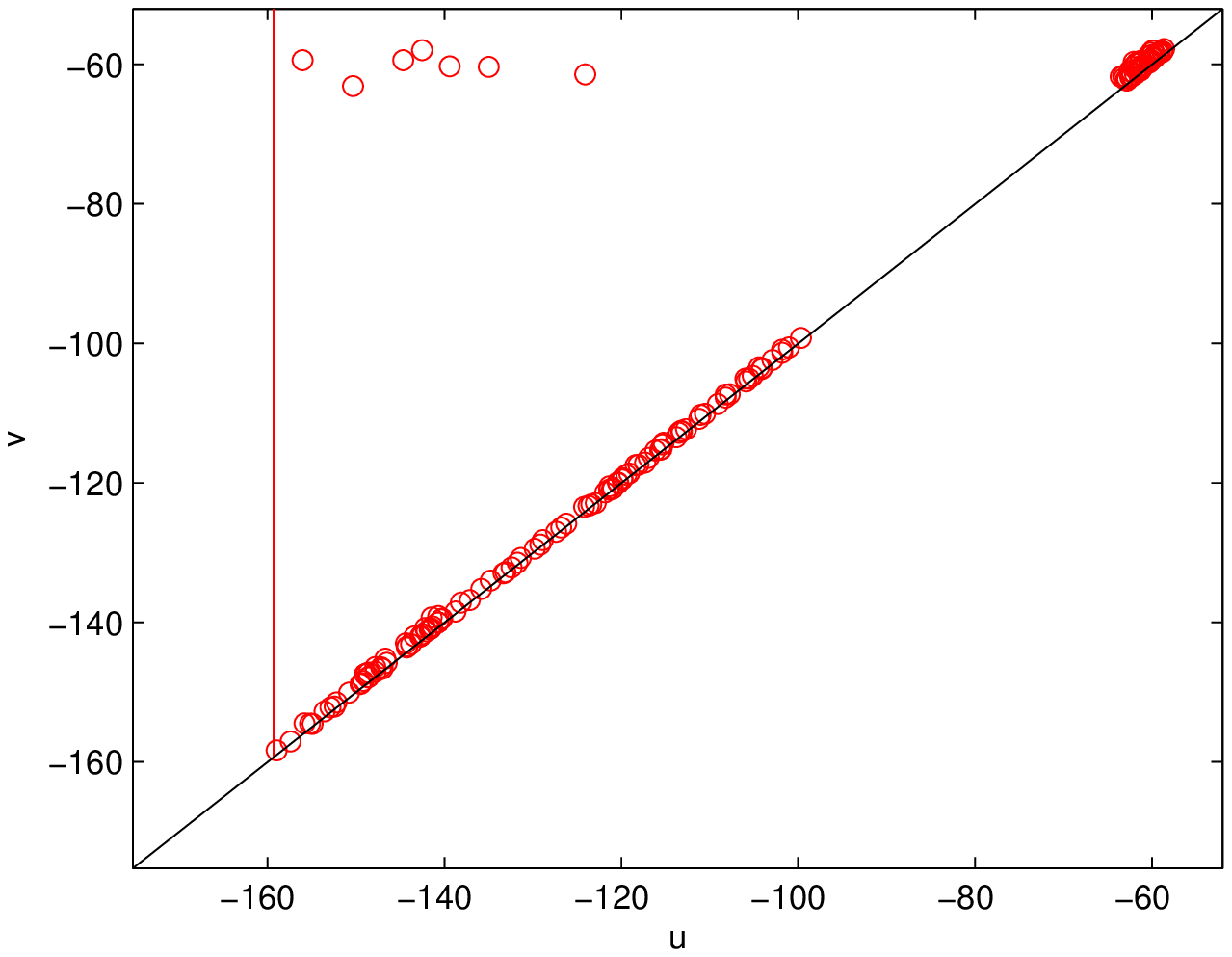} &  \includegraphics[width=0.5\textwidth]{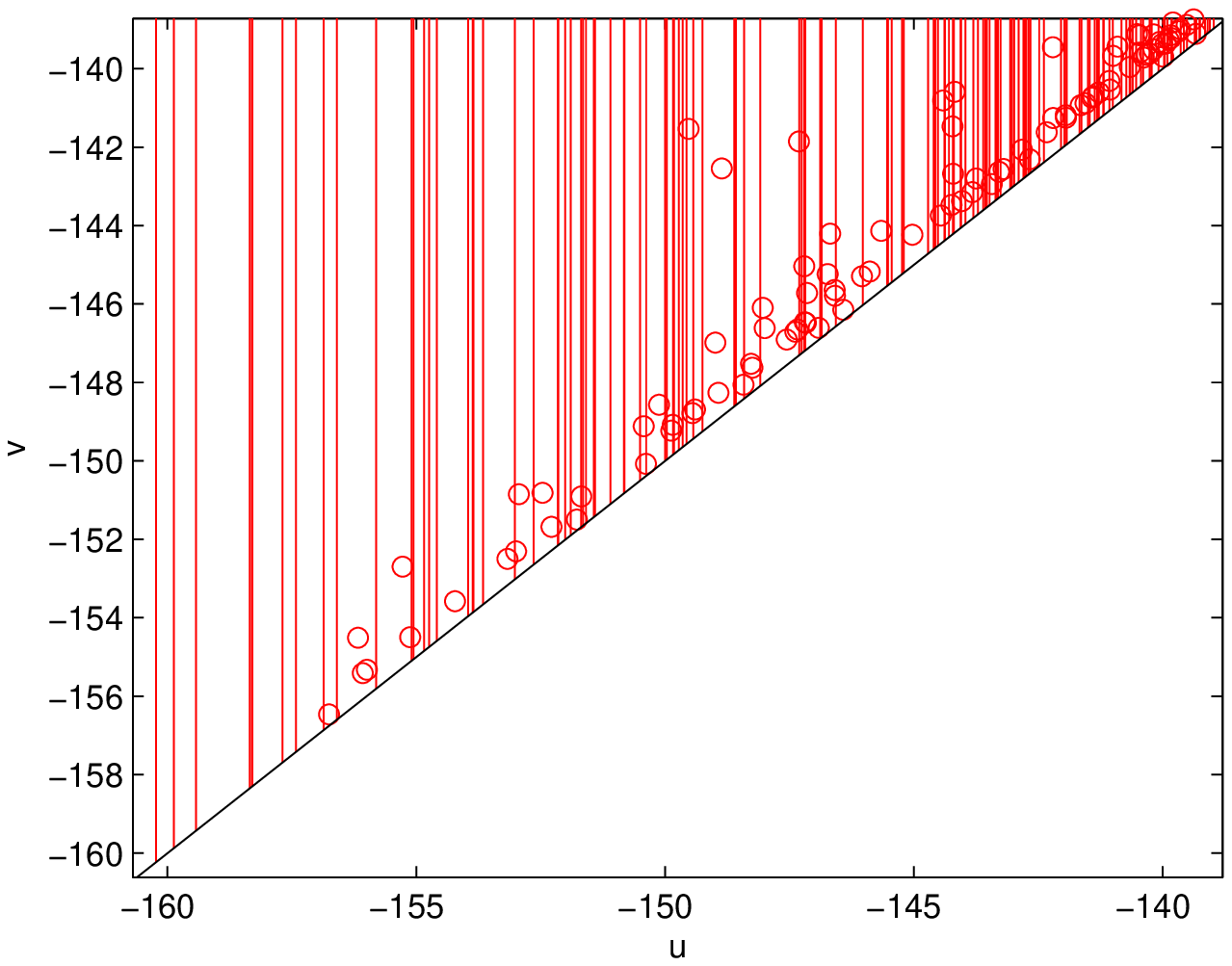}\end{tabular}
\caption{Left: The persistence diagram of the rank invariant $\rho_{(K_1,\p_{|K_1}),0}$ corresponding to the original octopus image. Right: A detail of the persistence diagram of the rank invariant $\rho_{(K_2,\p_{|K_2}),0}$ corresponding to the noisy octopus image.}\label{mono}
\end{figure}
\end{center}

Figure \ref{mono} (left) shows the persistence diagram of the $1$-dimensional $0$th rank invariant $\rho_{(K_1,\p_{|K_1}),0}$. It displays eight relevant points in the persistence diagram, corresponding to the eight tentacles of the octopus. Only one of these points is at infinity (and therefore depicted by a vertical line rather than by a circle) since $K_1$ has only one connected component. As for  $\rho_{(K_2,\p_{|K_2}),0}$, due to the presence of a great quantity of connected components in the noisy octopus,  its persistence diagram  has a very large number of points at infinity, and a figure showing them all would be hardly readable. For this reason Figure \ref{mono} (right) shows only a small subset of its  persistence diagram. However it is sufficient to perceive how dissimilar it is from $\rho_{(K_1,\p_{|K_1}),0}$.

\begin{center}
\begin{figure}[h]
\psfrag{s}{$s$}
\psfrag{t}{$t$}
\psfrag{ms}{$\mu\cdot s$}
\psfrag{mt}{$\mu\cdot t$}
\begin{tabular}{cc}
\includegraphics[width=0.5\textwidth]{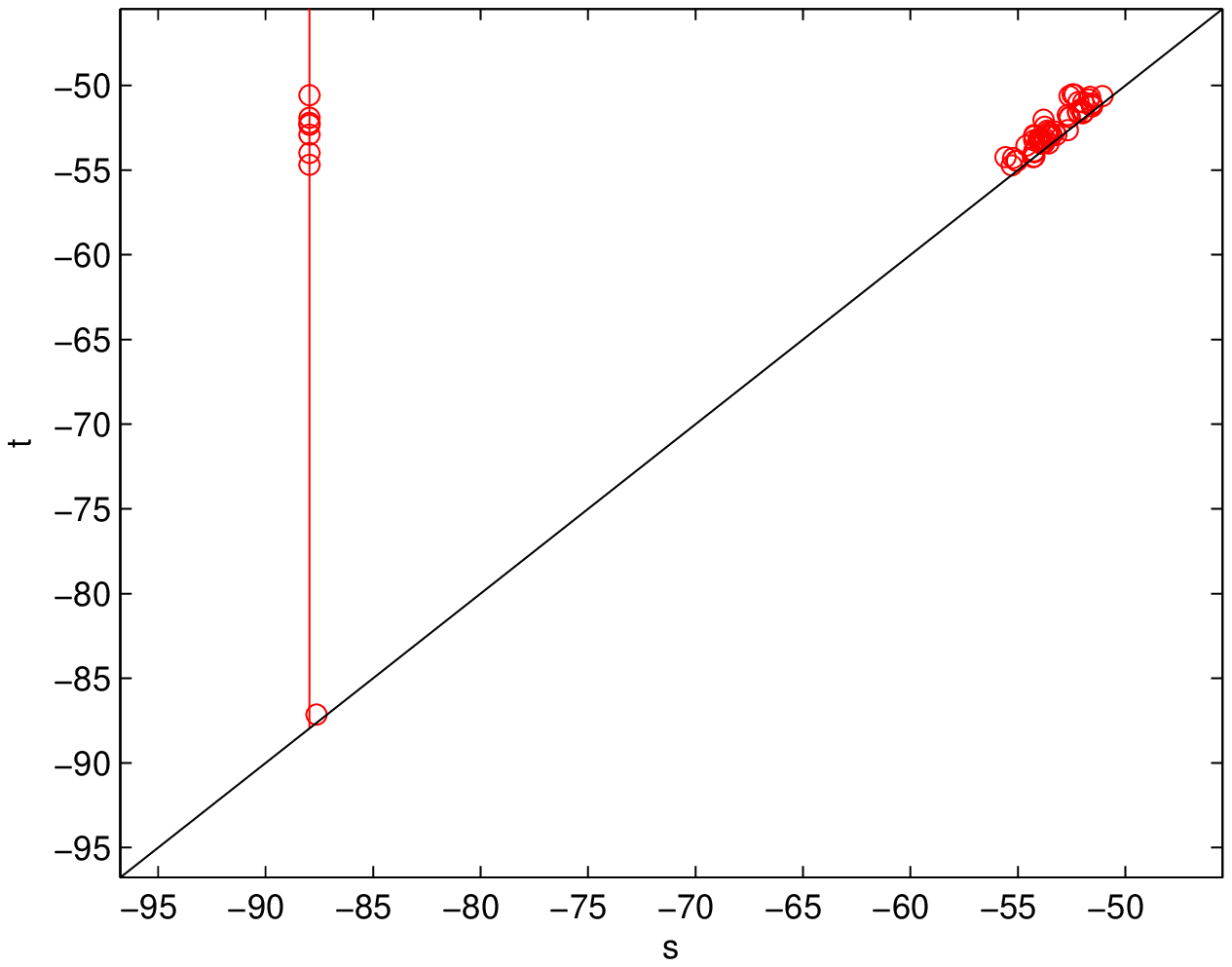} &  \includegraphics[width=0.5\textwidth]{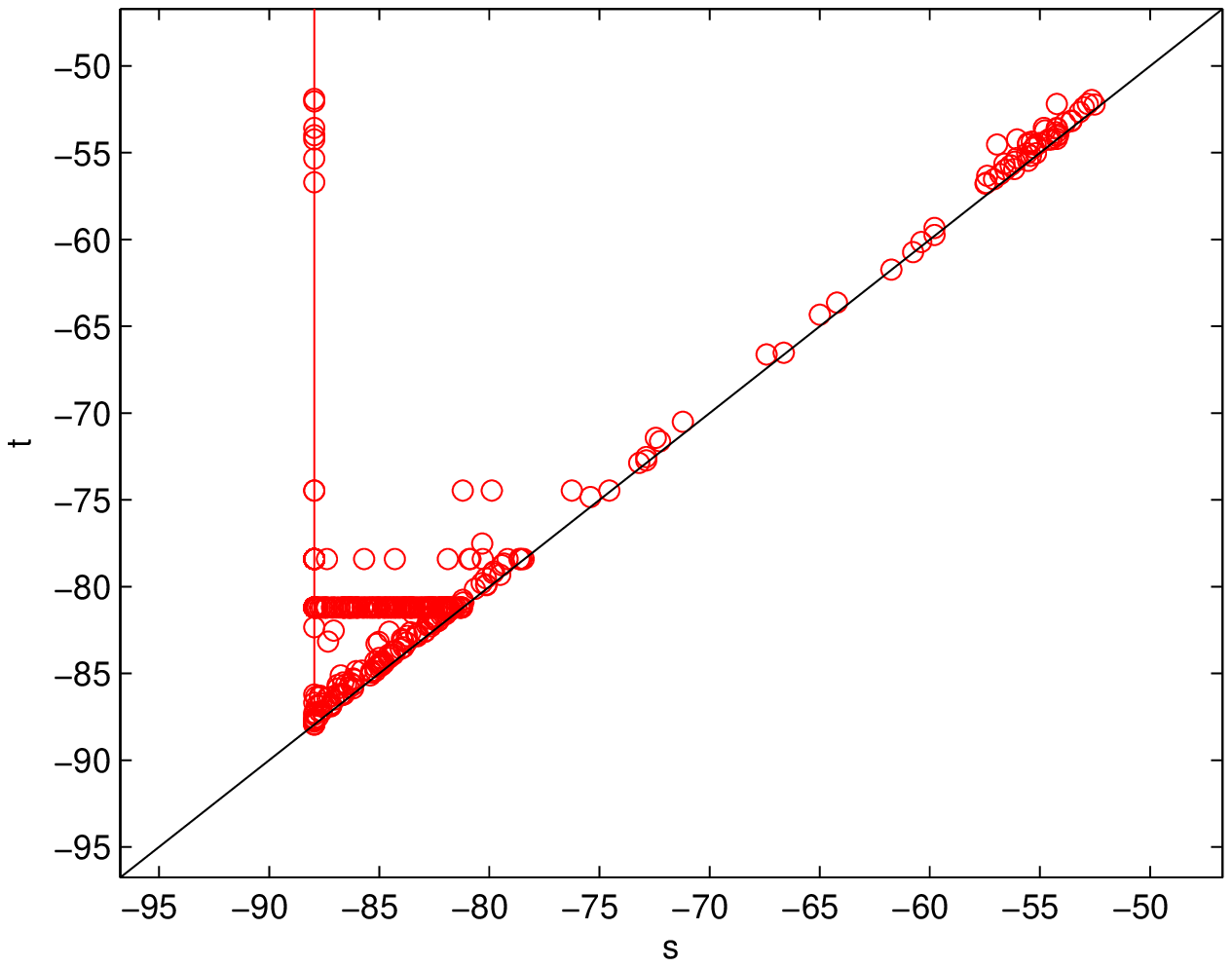}\\
(a) & (b)\\
\includegraphics[width=0.5\textwidth]{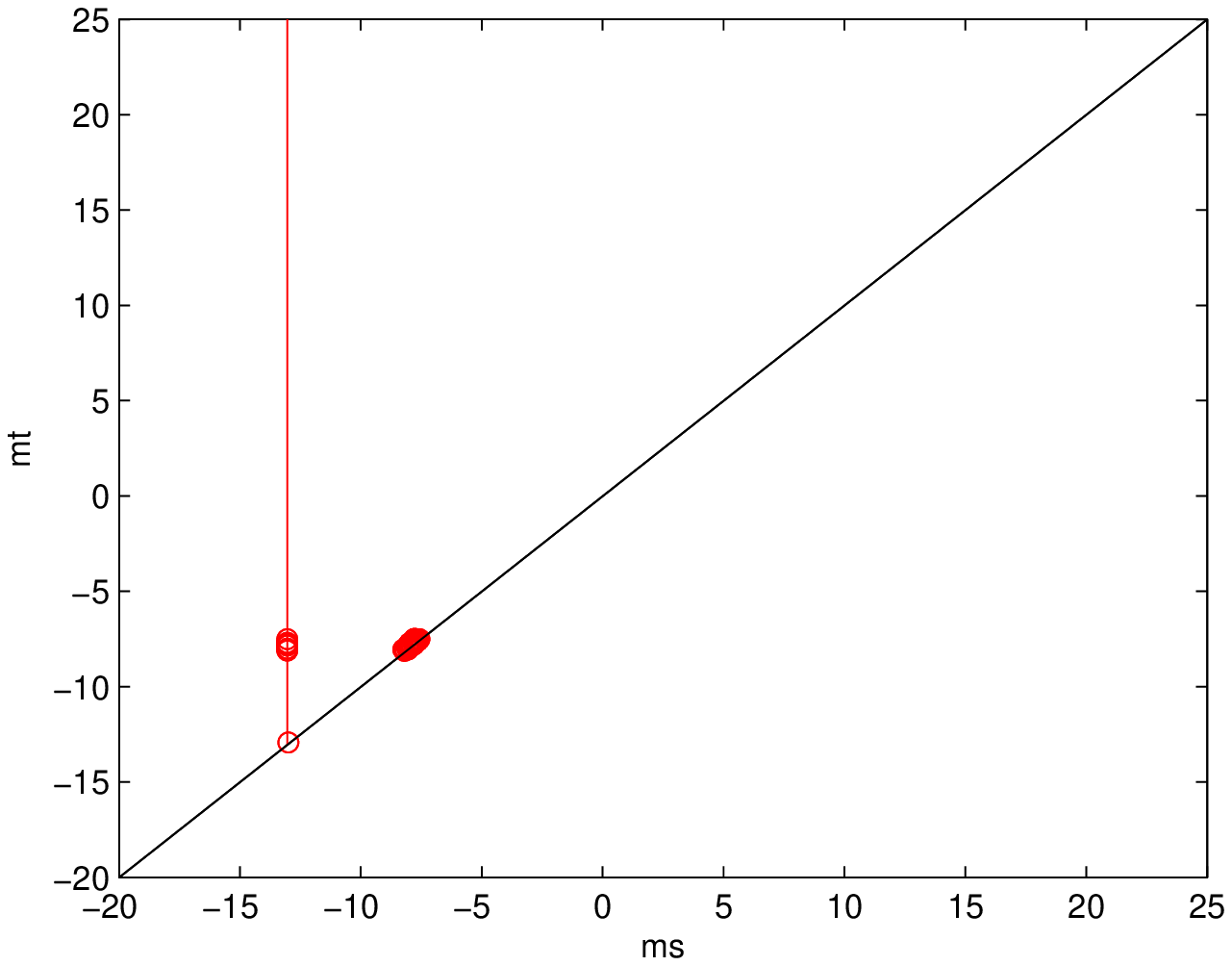} & \includegraphics[width=0.5\textwidth]{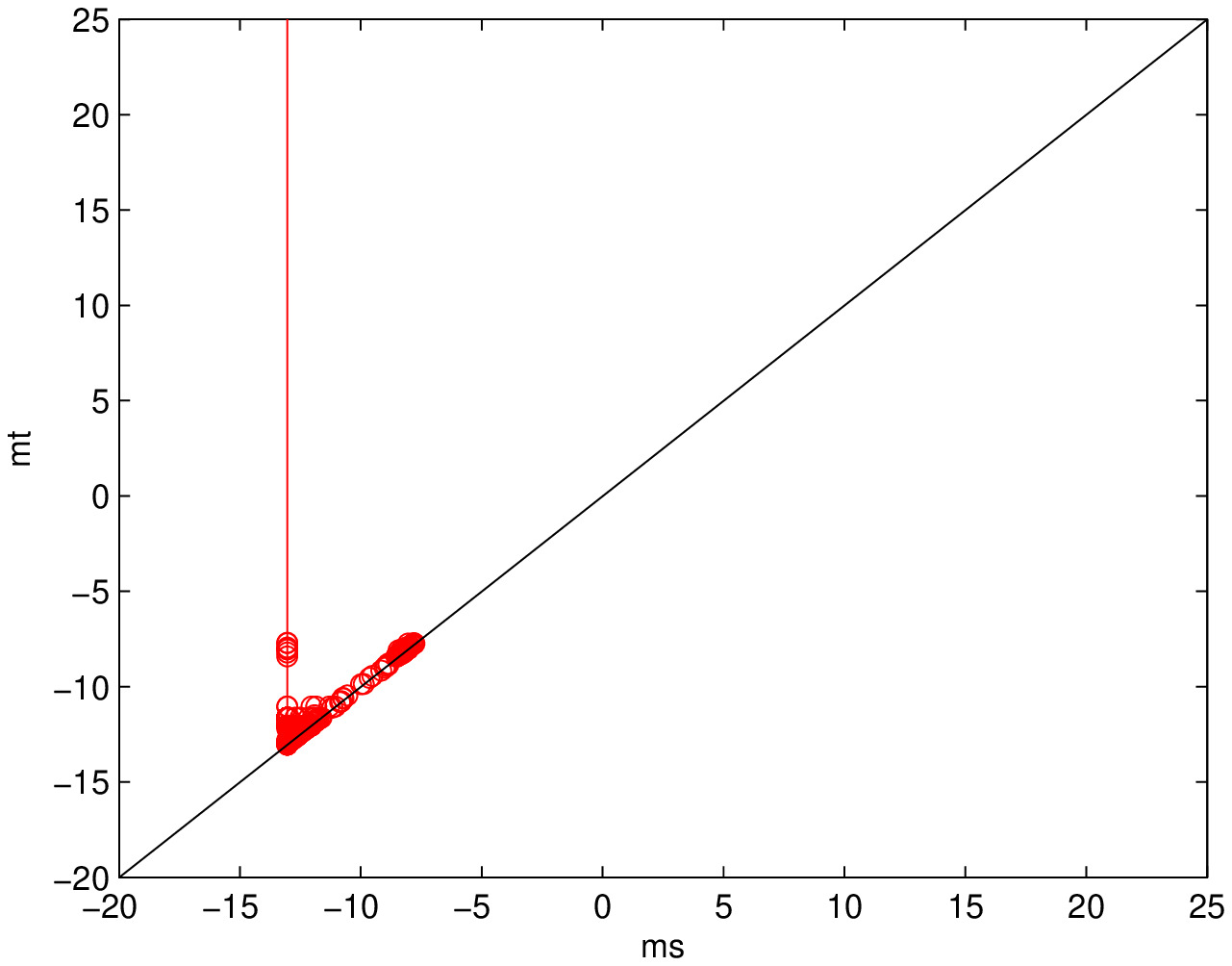}\\
 (c) & (d)\\
 \end{tabular}
\caption{(a) The rank invariant $\rho_{(D,\vec\Phi_1),0}$ restricted to the half-plane $\pi_{(\vec l,\vec b)}$,  with $\vec l=(0.1483, 0.9889)$ and $\vec b=(13.0434, -13.0434)$, that is the half-plane of the foliation containing the point $((0,-100),(3,-80))$. (b) The rank invariant $\rho_{(D,\vec\Phi_1),0}$ restricted to the same half-plane. (c)-(d) The same restrictions as in (a)-(b), respectively,  but rescaled by $\mu=\min\{l_1,l_2\}$. }\label{bidim1}
\end{figure}
\end{center}

As suggested by Theorem \ref{stability_domains}, if instead we compare $K_1$ and $K_2$ by means of the rank invariants  $\rho_{(D,\vec\Phi_1),0}$ and $\rho_{(D,\vec\Phi_2),0}$, where  $\vec \Phi_1:D\to \R^2$, $\vec \Phi_1=(d_{K_1},\p)$, and $\vec \Phi_2:D\to \R^2$, $\vec \Phi_2=(d_{K_2},\p)$, we can see the similarity between $K_1$ and $K_2$ modulo the {\em salt \& pepper} noise. This is illustrated in Figure \ref{bidim1}. In Figure \ref{bidim1} (a)-(b), we show the rank invariants $\rho_{(D,\vec\Phi_1),0}$ and $\rho_{(D,\vec\Phi_2),0}$ both restricted to the half-plane $\pi_{(\vec l,\vec b)}$,  with $\vec l=(0.1483, 0.9889)$ and $\vec b=(13.0434, -13.0434)$, that is the half-plane of the foliation containing the point $((0,-100),(3,-80))$.   In other words, Figure \ref{bidim1} (a)-(b) shows $\rho_{(D,F^{\vec \Phi_1}_{(\vec l,\vec b)}),0}$ and $\rho_{(D,F^{\vec \Phi_2}_{(\vec l,\vec b)}),0}$, respectively. We can appreciate their similarity, even if their matching distance $d_{match}$ is not necessarily smaller than the Hausdorff distance between $K_1$ and $K_2$. The considered half-plane has been chosen so that it contains points where the rank invariant takes non-trivial values.

Indeed it is easy to verify that Theorem \ref{stability_domains} does not guarantee the stability of $\rho_{(D,F_{(\vec
l,\vec b)}^{\vec \Phi}),q}$ but the stability of $\rho_{(D,\mu\cdot F_{(\vec
l,\vec b)}^{\vec \Phi}),q}$, where $\mu=\min_i l_i$. We point out that $\rho_{\left(D,\mu\cdot F_{(\vec
l,\vec b)}^{\vec \Phi}\right),q}\left(\mu\cdot s,\mu\cdot t\right)=\rho_{\left(D,F_{(\vec
l,\vec b)}^{\vec \Phi}\right),q}\left(s,t\right)$, and hence the passage from the measuring function $F_{(\vec
l,\vec b)}^{\vec \Phi}$ to the measuring function $\mu\cdot F_{(\vec
l,\vec b)}^{\vec \Phi}$ corresponds to ``rescaling up'' the domain of the rank invariant. In other words, when we change $K_1$ into a new compact set $K_2$ that is close to $K_1$ with respect to the Hausdorff distance, the matching distance between $\rho_{(K_1,\p_{|K_1}), q}$ and $\rho_{(K_2,\p_{|K_2}), q}$  may be not small, while the one between
$\rho_{\left(D,\mu\cdot F_{(\vec
l,\vec b)}^{(d_{K_1},\p_{|K_1})}\right), q}$ and $\rho_{\left(D,\mu\cdot F_{(\vec
l,\vec b)}^{(d_{K_2},\p_{|K_2})}\right), q}$ must be small.

This is illustrated in Figure \ref{bidim1}, where the rank invariants  $\rho_{(D, F^{\vec \Phi_1}_{(\vec l,\vec b)}),0}$ and $\rho_{(D, F^{\vec \Phi_2}_{(\vec l,\vec b)}),0}$, displayed the top row, are not as similar as the rank invariants $\rho_{(D,\mu\cdot F^{\vec \Phi_1}_{(\vec l,\vec b)}),0}$ and $\rho_{(D,\mu\cdot F^{\vec \Phi_2}_{(\vec l,\vec b)}),0}$, displayed in the bottom row.

Next we show how it is possible to point-wisely recover the rank invariant of $(K_1,\p_{|K_1})$ from that of $(D,\vec\Phi_1)$. According to  Theorem \ref{lim}, $\rho_{(K_1,\p_{|K_1}),0}(u,v)=\rho_{(D,\vec\Phi_1)}(\alpha,u,\beta,v)$ for $\alpha,\beta>0$ sufficiently small.

As shown in \cite{BiCe*08}, in this case the following equalities hold (with reference to Definition \ref{Admissible}):
\begin{eqnarray}\label{parameters}
\begin{array}{cc}
l_1(\alpha,u,\beta,v) = \frac{\beta-\alpha}{\sqrt{(\beta-\alpha)^2+(v-u)^2}}, & l_2(\alpha,u,\beta,v) = \frac{v-u}{\sqrt{(\beta-\alpha)^2+(v-u)^2}}, \\
b_1(\alpha,u,\beta,v) = \frac{\alpha(\beta+v)-\beta(\alpha +u)}{(\beta-\alpha)+(v-u)}, & b_2(\alpha,u,\beta,v) = \frac{u(\beta+v)-v(\alpha +u)}{(\beta-\alpha)+(v-u)},\\
s(a,u,\beta,v)   = \frac{\alpha +u}{l_1+l_2}, & t(\alpha,u,\beta,v) = \frac{\beta+v}{l_1+l_2}.
\end{array}
\end{eqnarray}

As a consequence, Theorems \ref{lim} and \ref{Reduction} (applied in this order) imply that for every pair $(u,v)$, with $u<v$, and for $0<\alpha<\beta$, with $\alpha$ and $\beta$ sufficiently small,  
\begin{eqnarray*}
  \rho_{(K_1,\p_{|K_1}), q}(u,v) & = &
  \rho_{(D,\vec\Phi_1),q}\left((\alpha,u),(\beta,v)\right)\\
   & = & \rho_{\left(D,F_{(\vec
l,\vec b)}^{\vec \Phi_1}\right),q}\left(s(\alpha,u,\beta,v),t(\alpha,u,\beta,v)\right)\\
 & = &
\rho_{\left(D,F_{(\vec
l,\vec b)}^{\vec \Phi_1}\right),q}\left(\frac{\alpha +u}{l_1+l_2},\frac{\beta+v}{l_1+l_2}\right)
\end{eqnarray*}
where $F_{(\vec
l,\vec b)}^{\vec \Phi_1}:D\rightarrow\R$ is defined by setting, for every $x\in D$,
$$
F_{(\vec
l,\vec b)}^{\vec \Phi_1}(x)=\max\left\{\frac{d_K(x)-b_1}{l_1},\frac{\varphi(x)-b_2}{l_2}\right\}\
.
$$

Hence the finite value $\rho_{(K_1,\p_{|K-1}), q}(u,v)$ is equal to $\rho_{\left(D,F_{(\vec
l,\vec b)}^{\vec \Phi_1}\right),q}\left(\frac{\alpha +u}{l_1+l_2},\frac{\beta+v}{l_1+l_2}\right)$, if we choose $\alpha$ and $\beta$ small enough. The corresponding admissible pair $(\vec l,\vec b)$ results to be close  to the pair $\left((0,1),(0,0)\right)$.

In other words, the information about the rank invariant of the original pair $(K_1,\p_{|K_1})$ can be recovered on the leaves associated with the admissible pairs $(\vec l,\vec b)$ in a small neighborhood of the pair $\left((0,1),(0,0)\right)$, after re-parameterizing these leaves by taking each point $(u,v)$ to the point $\left(\frac{\alpha +u}{l_1+l_2},\frac{\beta+v}{l_1+l_2}\right)$.
Note that the pair $\left((0,1),(0,0)\right)$ is not admissible but is located on the boundary of the set $Adm_2$. This leads to instabilities if we take $\alpha,\beta$ too small.

We underline that Theorem \ref{lim} ensures this approximation is  good only point-by-point. Thus, even for admissible pairs $(\vec l,\vec b)$ close enough to the pair $\left((0,1),(0,0)\right)$, the matching distance between $\rho_{\left(D,F_{(\vec
l,\vec b)}^{\vec \Phi_1}\right),q}$ and $\rho_{(K,\p_{|K_1}), q}$ may be quite large. \\

 To illustrate these issues we have computed the value taken by $\rho_{(K_1,\p_{|K_1}),0}$ at the point $(u,v)=(-100,-80)$, that is  $8$.  Using formulas (\ref{parameters}), we have computed the admissible pairs $\vec l=\vec l(\alpha,\beta)$, $\vec b=\vec b(\alpha,\beta)$ such that the half-plane $\pi_{(\vec l,\vec b)}$ contains the point $((\alpha,-100),(\beta,-80))$ for the values of $\alpha,\beta$ shown in Table \ref{table1}. For the same values of $\alpha,\beta$, Table \ref{table1} also shows the parameters $s=s(\alpha,\beta)$ and $t=t(\alpha,\beta)$
 for which we have that  $\rho_{(D,\vec\Phi_1),0}((\alpha,-100),(\beta,-80))=\rho_{\left(D,F_{(\vec
l,\vec b)}^{\vec \Phi_1}\right),0}\left(s,t\right)$.

Computations show that $\rho_{(K_1,\p_{|K_1}),0}(-100,-80)=\rho_{\left(D,F_{(\vec
l,\vec b)}^{\vec \Phi_1}\right),0}\left(s,t\right)=8$ for small but positive values of $\alpha$ and $\beta$. It is noticeable that for $\alpha=0$, due to the mentioned instabilities  near the boundary of the set $Adm_2$, computations are not reliable.

\begin{table}[h]
\begin{tabular}{|c|c|c|c|c|c|c|}\hline
$\alpha$ & $u$ & $\beta$ & $v$ &  $s$ & $t$ & $\rho_{\left(D,F_{(\vec
l(\eps),\vec b(\eps))}^{\vec \Phi_1}\right),0}\left(s,t\right)$\\
\hline
\vspace{-3mm} & & & & & & \\
0.5 & -100 & 24 & -80 & -70.5866 & -39.7272 &   1\\
0.5 & -100 & 16 & -80 & -70.9216 & -45.6179 &   3\\
0.5 & -100 & 8 & -80 & -77.2843 & -55.9262&   3\\
0.5 & -100 & 1 & -80 & -97.1120 & -77.1040 &   8\\
0.5 & -100 & 0.65 & -80 & -98.7692 & -78..7672 &   8\\
0.3 & -100 & 0.45 & -80 & -98.9677 & -78.9657 &   8\\
0.1 & -100 & 0.25 & -80 & -99.1663 & -79.1643 &   8\\
0 & -100 & 0.15 & -80 & -99.2655 & -79.2635 &   0\\
 \hline
\end{tabular}
\label{table1}
\vspace{5mm}\caption{The parameters used to approximate the value of $\rho_{(K_1,\p_{|K_1}),0}$ at $(u,v)=(-100,-80)$, that is $8$, using $\rho_{\left(D,F_{(\vec
l,\vec b)}^{\vec \Phi_1}\right),0}$.}
\end{table}

\begin{figure}[h]
\psfrag{s}{$s$}
\psfrag{t}{$t$}
\begin{tabular}{|c|c|}\hline
& \\
\begin{minipage}{.50\textwidth}
\centering
  $\alpha =0.5, \ \beta=24$\\
\includegraphics[width=.90\textwidth]{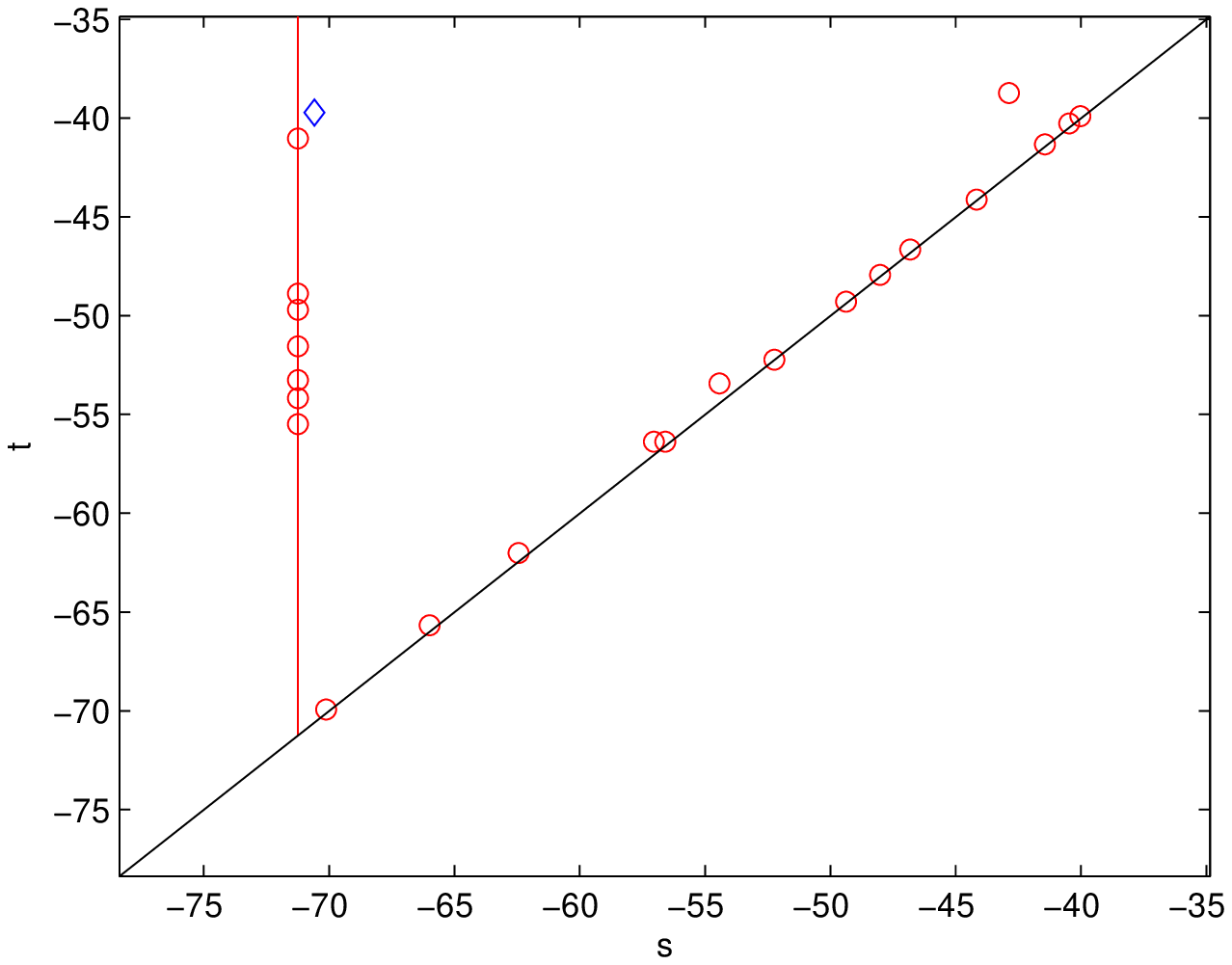}
\end{minipage}
&
\begin{minipage}{.50\textwidth}
\centering
 $\alpha =0.5, \ \beta=16$\\
\includegraphics[width=.90\textwidth]{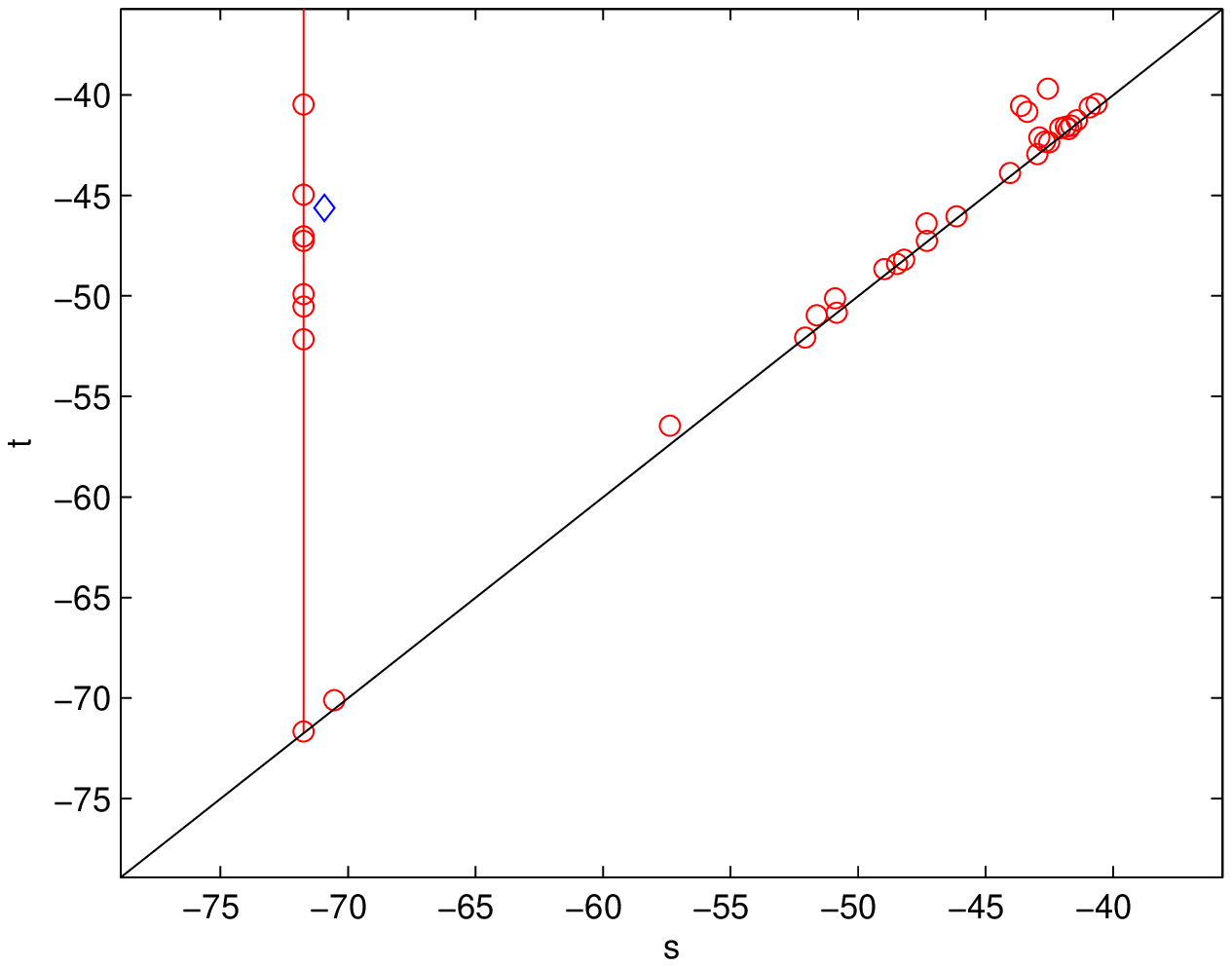}
\end{minipage}\\ \hline
& \\
\begin{minipage}{.50\textwidth}
\centering
$\alpha =0.5, \ \beta=8$\\
\includegraphics[width=.90\textwidth]{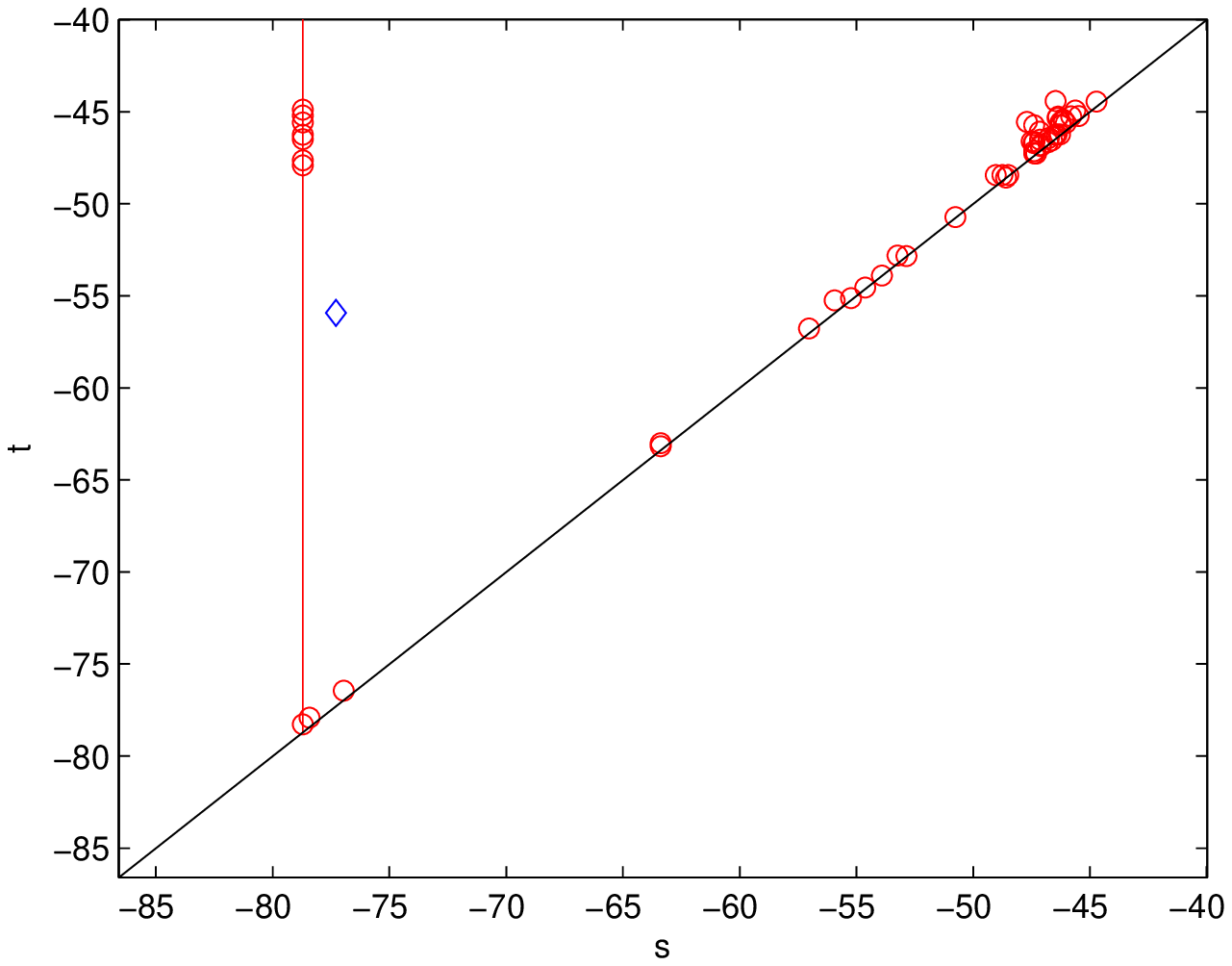}
\end{minipage}
&
\begin{minipage}{.50\textwidth}
\centering
$\alpha =0.5, \ \beta=1$\\
\includegraphics[width=.90\textwidth]{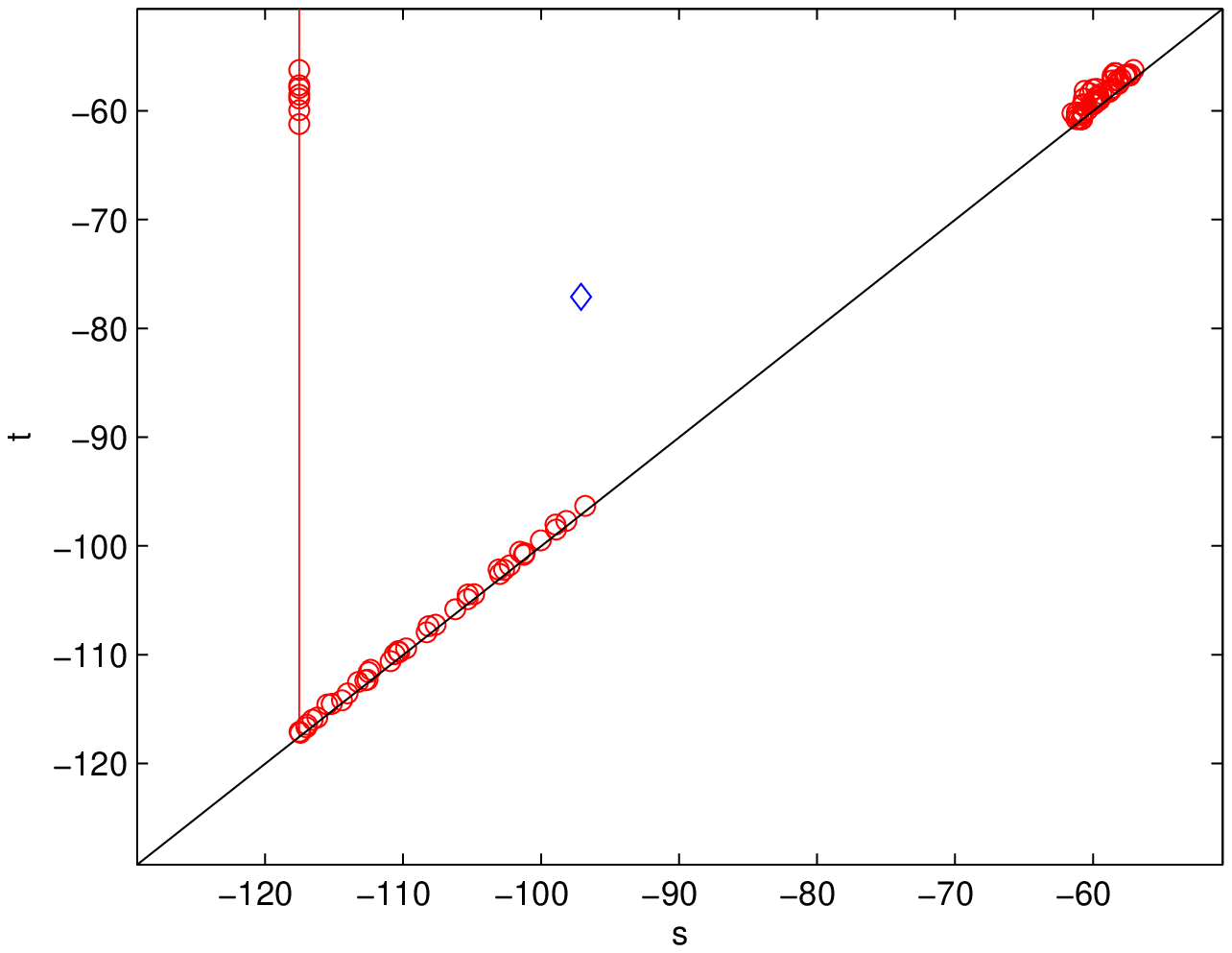}
\end{minipage}\\ \hline
& \\
\begin{minipage}{.50\textwidth}
\centering
$\alpha =0.5, \ \beta=0.65$\\
\includegraphics[width=.90\textwidth]{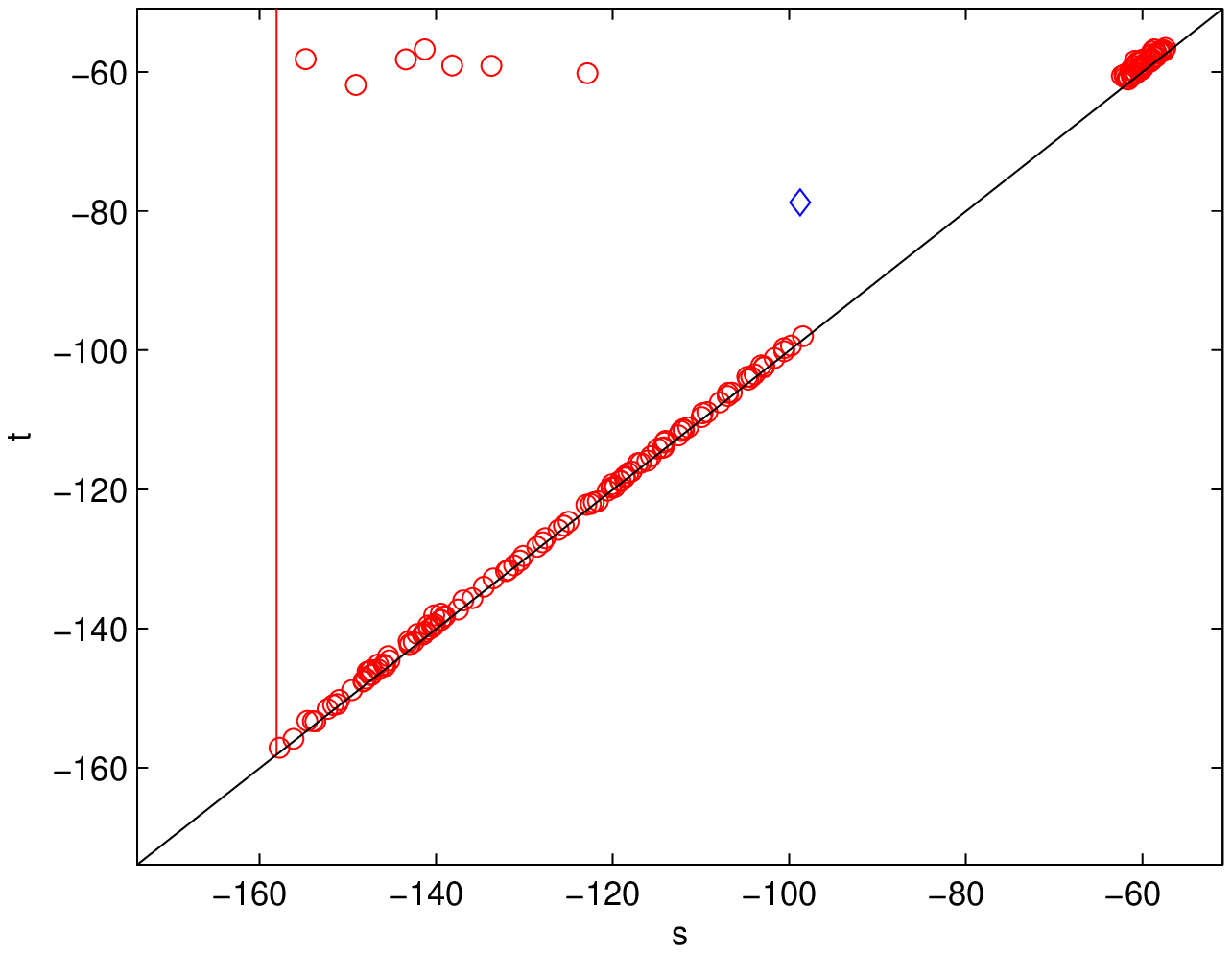}
\end{minipage}
&
\begin{minipage}{.50\textwidth}
\centering
$\alpha =0.3, \ \beta=0.45$\\
\includegraphics[width=.90\textwidth]{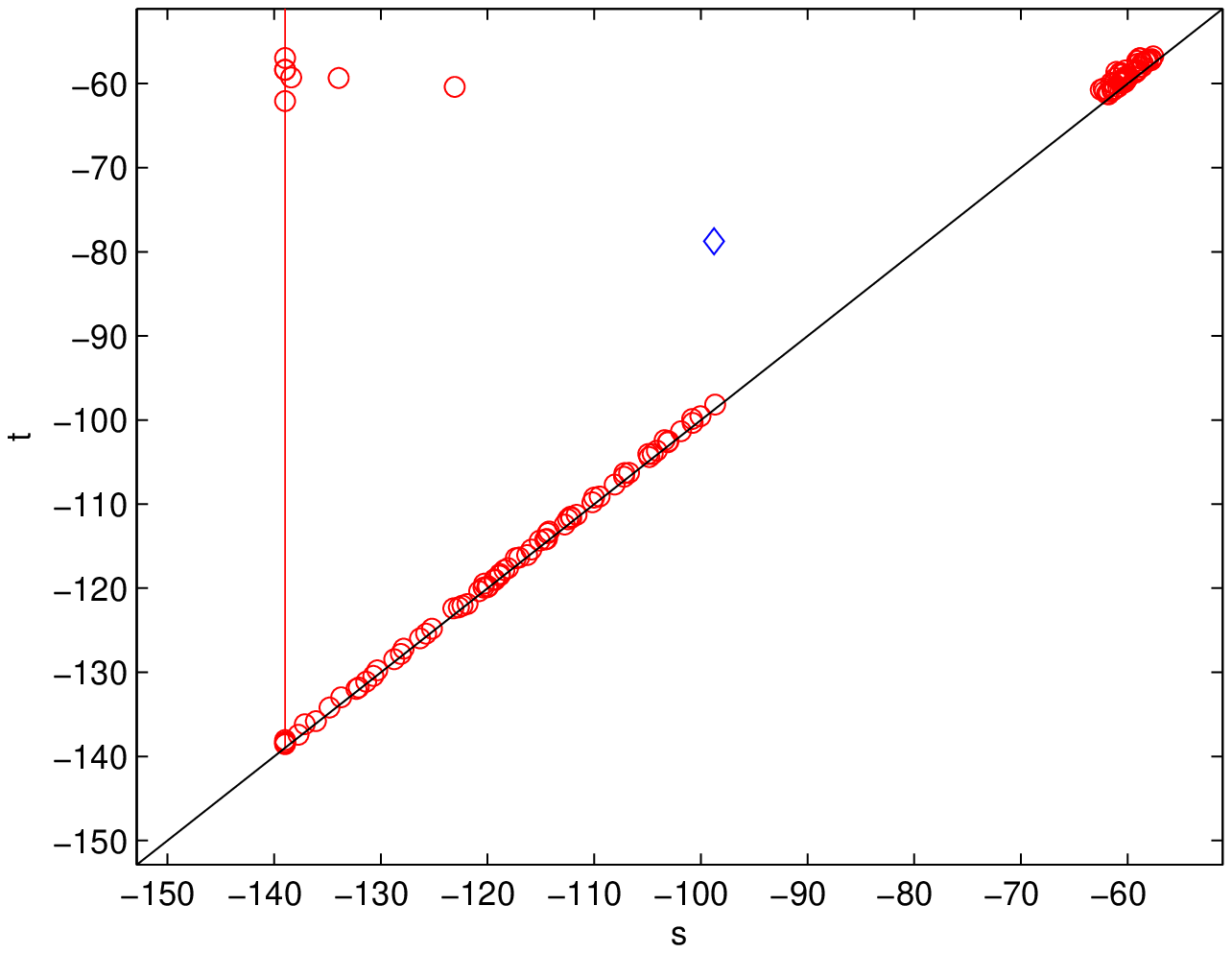}
\end{minipage}\\ \hline
& \\
\begin{minipage}{.50\textwidth}
\centering
$\alpha =0.1, \ \beta=0.25$\\
\includegraphics[width=.90\textwidth]{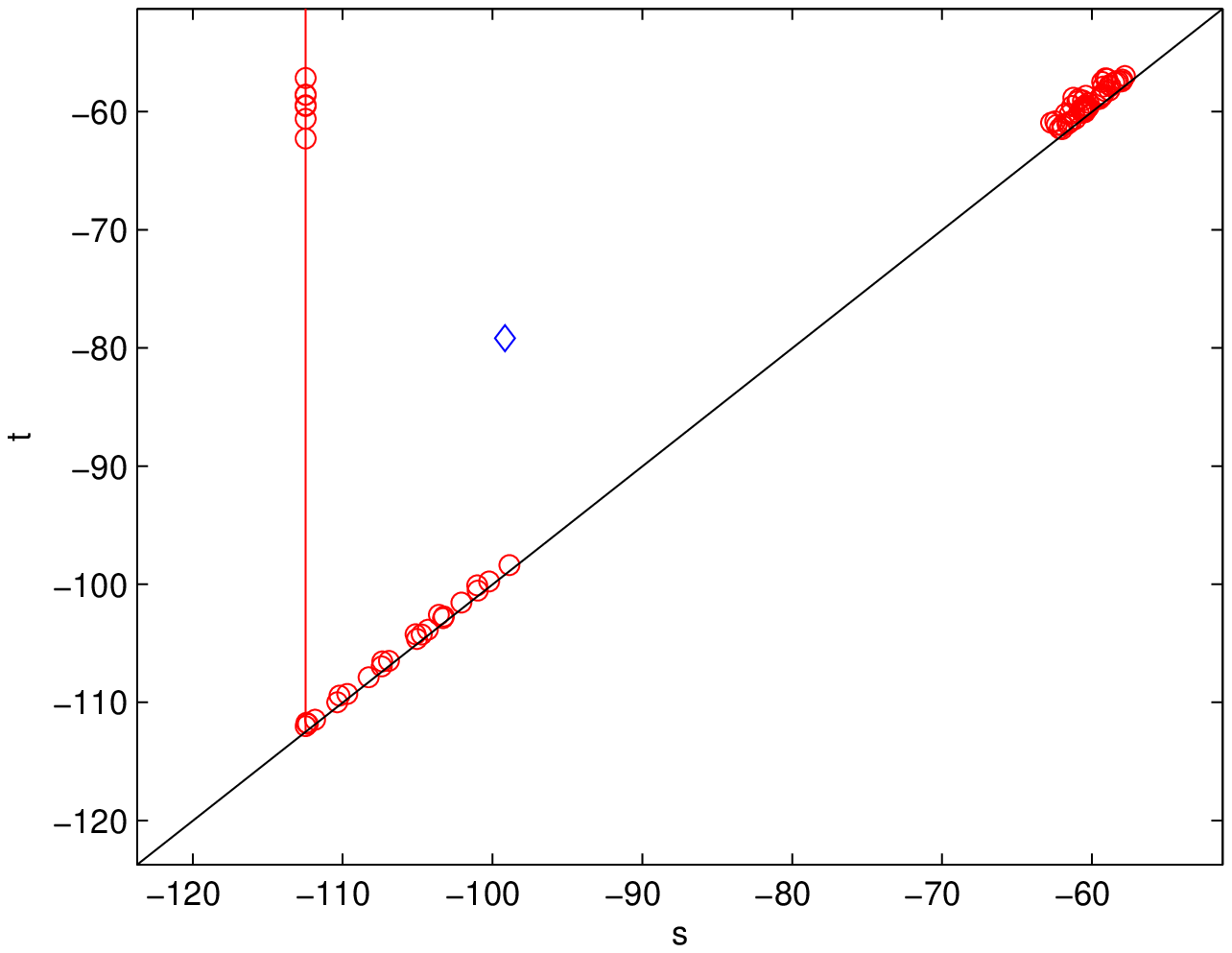}
\end{minipage}
&
\begin{minipage}{.50\textwidth}
\centering
$\alpha =0, \ \beta=0.15$\\
\includegraphics[width=.90\textwidth]{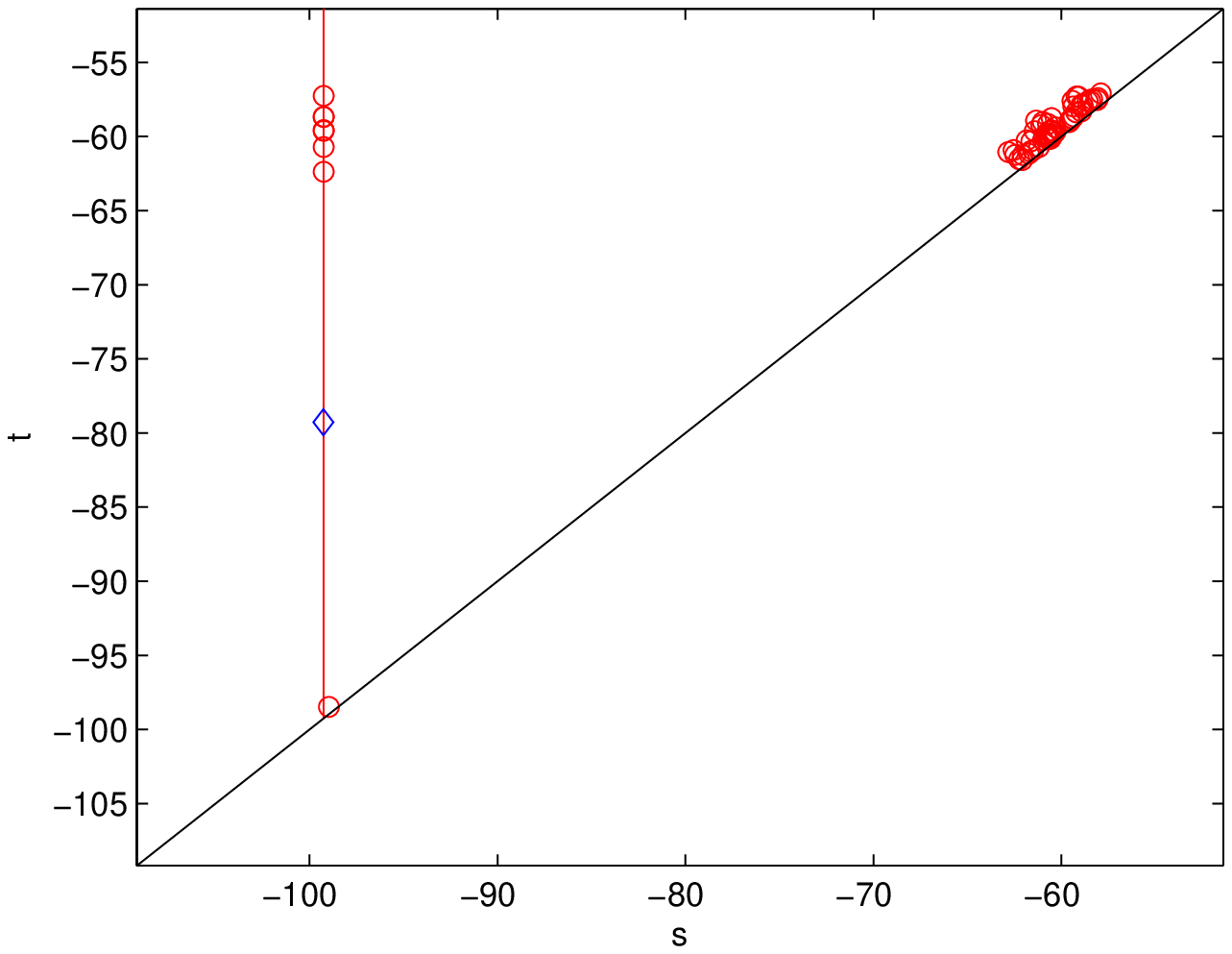}
\end{minipage}\\ \hline
\end{tabular}
\caption{  The rank invariant $\rho_{\left(D,F_{(\vec
l(\alpha),\vec b(\beta))}^{\vec \Phi_1}\right),0}$ as $\alpha$ and $\beta$ tend to $0$. Red circles and red lines  denote the points (proper or at infinity) of the corresponding persistence diagram. The blue diamonds denote the point $(s,t)$ corresponding to $((\alpha,-100),(\beta,-80))$.   }\label{bidim2}
\end{figure}

The corresponding rank invariants $\rho_{\left(D,F_{(\vec
l(\eps),\vec b(\eps))}^{\vec \Phi_1}\right),0}$ for the choices of $\alpha$ and $\beta$ considered in Table \ref{table1} are displayed in Figure \ref{bidim2}.

\section{Discussion}

In this paper we have shown the stability of persistent homology groups with respect to perturbations of the studied set. Measuring set perturbations by different distances requires different constructions in order to achieve stability. 

If set perturbations are measured through the Hausdorff distance, we replace compact sets  by distance functions. In this way, by the well-known property that if $K'$ is a good Hausdorff approximation of $K$ then the distance function $d_{K'}$ is close to $d_K$,   and by utilizing already available results of persistent homology stability  with respect to function perturbations, we deduce stability with respect to set perturbations. We also show that while passing from sets to functions we are still able  to recover information about the persistent homology groups of the original sets.  

If set perturbations are measured through the symmetric difference distance,  an analogous procedure leads to the proof of stability also in this case. Finally, using the $\sup$-distance  enables us to  guarantee stability also with respect to perturbations of fuzzy sets. 

The common underlying idea is to compare sets by comparing  functions describing the sets themselves.

While considering the Hausdorff distance and the symmetric difference distance is certainly not exhaustive of all the possible ways of measuring set perturbations, it accounts for two very widely used ones, and allows us to indicate a general procedure that could be applied also when dealing with other distances.

Finally, we underline that the technique developed in this paper essentially relies on the multidimensional generalization of persistent homology, showing once more the importance of further pursuing this area of research.

\bibliographystyle{abbrv}
\bibliography{biblio_stability_wrt_domain}
\end{document}